\documentclass[a4paper,11pt]{article}
\usepackage{latexsym,amssymb}
\usepackage{amsmath}
\usepackage[mathscr]{eucal}
\usepackage{color}
\usepackage[abbrev]{amsrefs}

\usepackage{latexsym}
\usepackage{amsmath}
\usepackage{amssymb}
\usepackage{enumerate}
\usepackage{mathrsfs}
\usepackage{color}

\newtheorem{Theorem}{\bf Theorem}[section]
\newtheorem{Lemma}{\bf Lemma}[section]
\newtheorem{Proposition}{\bf Proposition}[section]
\newtheorem{Corollary}{\bf Corollary}[section]
\newtheorem{Remark}{\bf Remark}[section]
\newtheorem{Example}{\bf Example}[section]
\newtheorem{Definition}{\bf Definition}[section]

\newenvironment{theorem}{\begin{Theorem}$\!\!\!$}{\end{Theorem}}
\newenvironment{lemma}{\begin{Lemma}$\!\!\!$}{\end{Lemma}}

\newenvironment{corollary}{\begin{Corollary}$\!\!\!$}{\end{Corollary}}

\newenvironment{definition}{\begin{Definition}$\!\!\!$}{\end{Definition}}

\numberwithin{equation}{section}

\allowdisplaybreaks[1]

\begin{document}

\title{Supercritical H\'enon type equation with a forcing term}
\author{
\qquad\\
Kazuhiro Ishige and Sho Katayama\vspace{5pt}\\
Graduate School of Mathematical Sciences, The University of Tokyo,\\ 
3-8-1 Komaba, Meguro-ku, Tokyo 153-8914, Japan\
\vspace{15pt}
}
\date{}
\maketitle
\begin{abstract}
This paper is concerned with the structure of solutions to the elliptic problem for an H\'enon type equation with a forcing term
\begin{equation}
\tag{$\mbox{P}_\kappa$}
-\Delta u=\alpha(x)u^p+\kappa\mu \quad\mbox{in}\quad{\mathbb R}^N,
\quad
u>0\quad\mbox{in}\quad{\mathbb R}^N,
\quad
u(x)\to 0\quad\mbox{as}\quad |x|\to\infty, 
\end{equation}
where $N\ge 3$, $p>1$, $\kappa>0$, $\alpha$ is a positive continuous function in ${\mathbb R}^N\setminus\{0\}$, 
and $\mu$ is a nonnegative Radon measure in ${\mathbb R}^N$. 
Under suitable assumptions on the exponent~$p$, the coefficient~$\alpha$, and the forcing term~$\mu$, 
we give a complete classification of the existence/nonexistence of solutions to problem~($\mbox{P}_\kappa$). 
\end{abstract}
\vspace{25pt}
\noindent E-mail Addresses:

\smallskip
\noindent 
{\tt ishige@ms.u-tokyo.ac.jp} (K. I.)\\
{\tt katayama-sho572@g.ecc.u-tokyo.ac.jp} (S. K.)\\
\vspace{20pt}
\vspace{3pt}
\newpage
\section{Introduction}
We are interested in the structure of solutions to 
the elliptic problem for an H\'enon type equation with a forcing term
\begin{equation}
\tag{$\mbox{P}_\kappa$}
\label{eq:P}
\left\{
\begin{array}{ll}
-\Delta u=\alpha(x)u^p+\kappa\mu & \quad\mbox{in}\quad{\mathbb R}^N,\vspace{3pt}\\
u>0 & \quad\mbox{in}\quad{\mathbb R}^N,\vspace{3pt}\\
u(x)\to 0 & \quad\mbox{as}\quad|x|\to\infty,
\end{array}
\right.
\end{equation}
where $N\ge 3$, $p>1$, $\kappa>0$, and $\mu$ is a nontrivial nonnegative Radon measure in ${\mathbb R}^N$. 
Here $\alpha$ is a positive continuous function in ${\mathbb R}^N\setminus\{0\}$ such that 
\begin{equation}
\label{eq:1.1}
\limsup_{|x|\to 0}|x|^{-a} \alpha(x)<\infty,
\qquad
\limsup_{|x|\to \infty}|x|^{-b} \alpha(x)<\infty, 
\end{equation}
for some $a>-2$ and $b\in{\mathbb R}$. 
Nonlinear elliptic equations with forcing terms in ${\mathbb R}^N$  
arise naturally in the study of stochastic processes. 
In particular, problem~\eqref{eq:P} with $\alpha\equiv 1$ in ${\mathbb R}^N$ 
appeared in establishing some limit theorems for super-Brownian motion. 
See e.g. \cites{Be, L, CFY, DGL} for a brief history and background of problem~\eqref{eq:P}.

In this paper, under suitable assumptions on the exponent~$p$ and the forcing term~$\mu$, 
we show the existence of a threshold $\kappa^*>0$ with the following properties:
\begin{itemize}
  \item[(A1)] 
  problem~\eqref{eq:P} possesses a minimal solution if $0<\kappa<\kappa^*$; 
  \item[(A2)] 
  problem~\eqref{eq:P} possesses a unique solution if $\kappa=\kappa^*$; 
  \item[(A3)]
  problem~\eqref{eq:P} possesses no solutions if $\kappa>\kappa^*$.
\end{itemize}
The existence/nonexistence and the behavior of solutions to problem~\eqref{eq:P} have been studied in many papers
(see e.g. \cites{Ba01, Ba02, BCP, BaL, BN, Be, CFY, CHZ, DGL, DLY, DY, HMP, LG, MWL} and references therein). 
There are some related results on properties~(A1) and (A3), 
however there are no available results on property~(A2) even in the case of $\alpha\equiv 1$ in ${\mathbb R}^N$. 

This paper is motivated by the paper~\cite{IOS01}, 
which treats the elliptic problem for the scalar field equation with a forcing term 
\begin{equation}
\tag{S}
\label{eq:S}
\left\{
\begin{array}{ll}
-\Delta u+u=u^p+\kappa\mu & \quad\mbox{in}\quad{\mathbb R}^N,\vspace{3pt}\\
u>0 & \quad\mbox{in}\quad{\mathbb R}^N,\vspace{3pt}\\
u(x)\to 0 & \quad\mbox{as}\quad|x|\to\infty,
\end{array}
\right.
\end{equation}
where $N\ge 2$, $p>1$, and $\mu$ is a nontrivial nonnegative Radon measure in ${\mathbb R}^N$ with compact support. 
In \cite{IOS01}, 
under the following condition on $\mu$:
\begin{itemize}
\item
$G*\mu\in L^q({\mathbb R}^N)$ for some $q\in(p,\infty]$ with $q>N(p-1)/2$, 
where $G$ is the fundamental solution to the elliptic operator $-\Delta+1$ in ${\mathbb R}^N$,
\end{itemize}
the existence of a threshold with properties~(A1)--(A3) was proved for problem~\eqref{eq:S} in the case of $1<p<p_{JL}$,
where $p_{JL}$ is the Joseph--Lundgren exponent, that is, 
\begin{equation}
\label{eq:1.2}
p_{JL}:=
\left\{
\begin{array}{ll}
1+\displaystyle{\frac{4}{N-4-\sqrt{4N-4}}} & \mbox{if}\quad N>10,\vspace{5pt}\\
\infty & \mbox{otherwise}.
\end{array}
\right.
\end{equation}
Unfortunately, the arguments in \cite{IOS01} depend heavily on the exponential decay of the fundamental solution $G$ at the space infinity, 
and they are not applicable to our problem~\eqref{eq:P}. 
See e.g. \cites{IOS02, IOS03, NS01, NS02} for related results on problem~(S). 
\medskip

We introduce some notation. 
In what follows, unless otherwise stated, let $N\ge 3$. 
For any $R>0$, let 
$$
B(0,R):=\{y\in{\mathbb R}^N\,:\,|y|<R\},\quad 
A(0,R):=B(0,R)\setminus B(0,R/2).
$$
We denote by ${\mathcal M}_+$ the set of all nontrivial nonnegative Radon measures in ${\mathbb R}^N$. 
Let $\Gamma$ be the fundamental solution to $-\Delta v=0$  in ${\mathbb R}^N$, that is,
$$
\Gamma(x):=\frac{1}{N(N-2)\omega_N}|x|^{-N+2},
$$
where $\omega_N$ is the volume of the ball $B(0,1)$. 
We denote by ${\mathcal D}^{1,2}$ the completion of $C_c^\infty({\mathbb R}^N)$ with respect to the norm $\|\nabla\cdot\|_{L^2({\mathbb R}^N)}$. 
It follows from the Sobolev inequality that 
$$
{\mathcal D}^{1,2}
=\left\{v\in L^\frac{2N}{N-2}({\mathbb R}^N)\cap W^{1,2}_\mathrm{loc}({\mathbb R}^N)\colon\int_{{\mathbb R}^N}|\nabla v|^2\,dx<\infty\right\}.
$$
For any $q\in[1,\infty]$ and $\beta$, $\gamma\in{\mathbb R}$, 
we define
$$
L^q_{\beta,\gamma}:=\left\{f\in L^q_{{\rm loc}}({\mathbb R}^N\setminus\{0\})\,:\,\|f\|_{L^q_{\beta,\gamma}}<\infty\right\},
$$
where
$$
\|f\|_{L^q_{\beta,\gamma}}
:=\sup_{R>0}R^{-\frac{N}{q}}\frac{\|f\|_{L^q(A(0,R))}}{\omega_{\beta,\gamma}(R)},
\qquad
\omega_{\beta,\gamma}(R):=
\left\{
\begin{array}{ll}
R^\beta & \mbox{for $R\in(0,1]$},\vspace{5pt}\\
R^\gamma & \mbox{for $R\in(1,\infty)$}.
\end{array}
\right.
$$
Notice that, if $f\in C({\mathbb R}^N\setminus\{0\})$ satisfies $f(x)=O(|x|^\beta)$ as $|x|\to 0$ 
and $f(x)=O(|x|^\gamma)$ as $|x|\to \infty$, then $f\in L^q_{\beta,\gamma}$ for any $q\in[1,\infty]$. 
\vspace{5pt}

We define solutions to problem~\eqref{eq:P}.
\begin{definition}
\label{Definition:1.1} 
Let $\mu\in{\mathcal M}_+$ and $\kappa>0$.  
\begin{enumerate}
\item[{\rm (1)}] 
Let $u$ be a nonnegative, measurable, finite and positive almost everywhere function in ${\mathbb R}^N$. 
We say that $u$ is a solution to problem~{\rm \eqref{eq:P}} if 
$u$~satisfies 
$$
u(x)=\int_{{\mathbb R}^N}\Gamma(x-y)\alpha(y)u(y)^p\,dy+\kappa \int_{{\mathbb R}^N}\Gamma(x-y)\,d\mu(y)
$$
for almost all {\rm(}a.a.{\rm)}~$x\in{\mathbb R}^N$.
We also say that $u$ is a supersolution to problem~{\rm \eqref{eq:P}} 
if~$u$ satisfies 
\[
u(x)\ge \int_{{\mathbb R}^N}\Gamma(x-y)\alpha(y)u(y)^p\,dy+\kappa \int_{{\mathbb R}^N}\Gamma(x-y)\,d\mu(y) 
\]
for a.a.~$x\in{\mathbb R}^N$.
\item[{\rm (2)}]
Let $u$ be a solution to problem~{\rm \eqref{eq:P}}. 
We say that $u$ is a minimal solution to problem~{\rm \eqref{eq:P}} 
if, for any solution $v$ to problem~{\rm \eqref{eq:P}}, the inequality $u(x)\le v(x)$ holds for a.a.~$x\in{\mathbb R}^N$.  
{\rm ({\it Obviously, a minimal solution is uniquely determined.})}
\end{enumerate}
\end{definition}

Now we are ready to state our results on problem~\eqref{eq:P}. 
Theorem~\ref{Theorem:1.1} concerns properties~(A1) and (A3) and 
it shows the existence of a threshold 
on the existence/nonexistence of solutions to problem~\eqref{eq:P}.
\begin{theorem}
\label{Theorem:1.1}
Let $N\ge 3$. Let $\alpha$ be a positive continuous function in ${\mathbb R}^N\setminus\{0\}$ and assume~\eqref{eq:1.1}. 
Let $1<p<r\le\infty$ and $c$, $d\in{\mathbb R}$ be such that 
\begin{equation}
\label{eq:1.3}
p>\frac{N+b}{N-2},\qquad r>\frac{N(p-1)}{2},
\qquad c>-\frac{a+2}{p-1},
\qquad d<-\frac{b+2}{p-1}.
\end{equation}
Assume that $\mu\in{\mathcal M}_+$ satisfies
\begin{equation}
\label{eq:1.4}
\Gamma*\mu\in L^r_{c,d}.
\end{equation}
Then there exists $\kappa^*>0$ such that
\begin{itemize}
  \item[{\rm (1)}] 
  if $0<\kappa<\kappa^*$, then problem~\eqref{eq:P} possesses a minimal solution $u^\kappa$;
  \item[{\rm (2)}] 
  if $\kappa>\kappa^*$, then problem~\eqref{eq:P} possesses no solutions.
\end{itemize}
Furthermore, there exists $\kappa_*\in(0,\kappa^*]$ such that
\begin{itemize}
 \item[{\rm (3)}] 
if $0<\kappa<\kappa_*$, then 
$u^\kappa\in L^r_{c_*,d_*}$ with $c_*:=\min\{c,0\}$ and $d_*:=\max\{d,-N+2\}$. 
\end{itemize}
\end{theorem}
We notice that, 
under condition~\eqref{eq:1.1}, 
problem~\eqref{eq:P} possesses no solutions generally if 
$$
1<p\le\frac{N+b}{N-2}.
$$
See e.g. \cite{BP}*{Theorem~3.3}. (See also \cite{Gidas} and \cite{QS}*{Section~8.1} for the case of $\alpha\equiv 1$ in ${\mathbb R}^N$.) 

Theorem~\ref{Theorem:1.2} concerns property~(A2) and it is the main result of this paper. 
Let
\begin{align*}
p^*(\eta) & :=
\left\{
\begin{array}{ll}
1+\displaystyle{\frac{2(\eta+2)}{N-\eta-4-\sqrt{(\eta+2)(2N+\eta-2)}}} & \mbox{if}\quad N>10+4\eta,\vspace{5pt}\\
\infty & \mbox{otherwise},
\end{array}
\right.\\
p_*(\eta) & :=\left\{
\begin{array}{ll}
1+\displaystyle{\frac{2(\eta+2)}{N-\eta-4+\sqrt{(\eta+2)(2N+\eta-2)}}} & \mbox{if}\quad \eta>-2,\vspace{5pt}\\
1 & \mbox{otherwise},
\end{array}
\right.
\end{align*}
for $\eta\in{\mathbb R}$. 
Then $p_{JL}=p^*(0)$ and 
\begin{equation}
\label{eq:1.5}
\begin{split}
 & N-\eta-4-\sqrt{(\eta+2)(2N+\eta-2)}>0\quad\mbox{if $N>10+4\eta$},\\
 & N-\eta-4+\sqrt{(\eta+2)(2N+\eta-2)}>0\quad\mbox{if $\eta>-2$}.
\end{split}
\end{equation}
\begin{theorem}
\label{Theorem:1.2}
Assume the same conditions as in Theorem~{\rm\ref{Theorem:1.1}}. 
Let $\kappa_*$ and $\kappa^*$ be as in Theorem~{\rm\ref{Theorem:1.1}}. 
Further, assume that 
$$
p_*(b)<p<p^*(a_-)\quad\mbox{with}\quad a_-:=\min\{a,0\}.
$$
Then $\kappa^*=\kappa_*$ and the following properties hold.
\begin{itemize}
  \item[{\rm (1)}] 
  If $0<\kappa\le\kappa^*$, then problem~\eqref{eq:P} possesses a minimal solution $u^\kappa\in L^q_{c_*,d_*}$. 
  Furthermore, if $\kappa=\kappa^*$, then $u^{\kappa^*}$ is a unique solution to problem~\eqref{eq:P}.
  \item[{\rm (2)}] 
  If $\kappa>\kappa^*$, then problem~\eqref{eq:P} possesses no solutions.
\end{itemize}
\end{theorem}

As corollaries of Theorems~\ref{Theorem:1.1} and \ref{Theorem:1.2}, we have
the following results in the case of $\alpha\equiv 1$ in~${\mathbb R}^N$, that is, 
the elliptic problem for the Lane-Emden equation with a forcing term, 
\begin{equation}
\tag{$\mbox{P}'_\kappa$}
\label{eq:P'}
\left\{
\begin{array}{ll}
-\Delta u=u^p+\kappa\mu & \quad\mbox{in}\quad{\mathbb R}^N,\vspace{3pt}\\
u>0 & \quad\mbox{in}\quad{\mathbb R}^N,\vspace{3pt}\\
u(x)\to 0 & \quad\mbox{as}\quad|x|\to\infty,
\end{array}
\right.
\end{equation}
where $N\ge 3$, $p>1$, $\kappa>0$, and $\mu\in{\mathcal M}_+$. 
\begin{corollary}
\label{Corollary:1.1}
Let $N\ge 3$. 
Let $1<p<r\le\infty$ be such that 
$$
p>\frac{N}{N-2},\qquad r>\frac{N(p-1)}{2}.
$$
Assume that $\mu\in{\mathcal M}_+$ satisfies
\begin{equation}
\label{eq:1.6}
\Gamma*\mu\in L^r_{{\rm loc}}({\mathbb R}^N),
\quad
\limsup_{R\to\infty}R^{-\frac{N}{r}-\theta}\|\Gamma*\mu\|_{L^r(A(0,R))}<\infty
\quad\mbox{for some}\quad \theta<-\frac{2}{p-1}.
\end{equation}
Then there exists $\kappa^*>0$ such that
\begin{itemize}
  \item[{\rm (1)}] 
  if $0<\kappa<\kappa^*$, then problem~\eqref{eq:P'} possesses a minimal solution $u^\kappa$;
  \item[{\rm (2)}] 
  if $\kappa>\kappa^*$, then problem~\eqref{eq:P'} possesses no solutions.
\end{itemize}
Furthermore, there exists $\kappa_*\in(0,\kappa^*]$ such that
\begin{itemize}
 \item[{\rm (3)}] 
if $0<\kappa<\kappa_*$, then 
\begin{equation}
\label{eq:1.7}
u^\kappa\in L^r_{{\rm loc}}({\mathbb R}^N),
\qquad
\limsup_{R\to\infty}R^{-\frac{N}{r}-\theta_*}\|u^\kappa\|_{L^r(A(0,R))}<\infty,
\end{equation}
where $\theta_*:=\max\{\theta,-N+2\}$. 
\end{itemize}
\end{corollary}
\begin{corollary}
\label{Corollary:1.2}
Assume the same conditions as in Theorem~{\rm\ref{Theorem:1.1}}. 
Further, assume that 
$$
p_{JL}'<p<p_{JL},
$$
where $p_{JL}$ is as in \eqref{eq:1.2} and 
$$
p_{JL}'
:=p_*(0)=1+\displaystyle{\frac{4}{N-4+\sqrt{4N-4}}}.
$$
Then $\kappa^*=\kappa_*$, where $\kappa^*$ and $\kappa_*$ are as in Corollary~{\rm\ref{Corollary:1.1}}. 
Furthermore, if $\kappa=\kappa^*$, then problem~{\rm (P')} possesses a unique solution $u^{\kappa^*}$ 
and it satisfies \eqref{eq:1.7} with $\kappa=\kappa^*$. 
\end{corollary}
Assumption~\eqref{eq:1.6} is somewhat more general than that of \eqref{eq:1.4} with $c=0$. 
\vspace{5pt}

In the proofs of our theorems, for any $\kappa>0$, 
we set $U^\kappa_{-1}\equiv 0$ in ${\mathbb R}^N$ and 
define approximate solutions $\{U^\kappa_j\}$ to problem~\eqref{eq:P} and their differences $\{V_j^\kappa\}$ inductively by 
\begin{equation}
\label{eq:1.8}
\begin{array}{ll}
U_j^\kappa(x):=[\Gamma*\alpha(U_{j-1}^\kappa)^p](x)+\kappa (\Gamma*\mu)(x), & \quad j=0,1,2,\dots,\vspace{5pt}\\
V_j^\kappa(x):=U_j^\kappa(x)-U_{j-1}^\kappa(x), & \quad j=0,1,2,\dots,
\end{array}
\end{equation}
for a.a.~$x\in{\mathbb R}^N$. 
By induction we easily see that 
\begin{equation}
\label{eq:1.9}
U_{j+1}^\kappa(x)>U_j^\kappa(x),\qquad
0<V_j^\kappa(x)\le p\left[\Gamma*\left(\alpha (U^\kappa_{j-1})^{p-1}V^\kappa_{j-1}\right)\right](x),
\end{equation}
for a.a.~$x\in{\mathbb R}^N$ and $j=0,1,2,\dots$. 
Assume the same conditions as in Theorem~\ref{Theorem:1.1}. 
Define
\begin{equation}
\label{eq:1.10}
\begin{split}
 & {\mathcal K}^*:=\left\{\kappa>0\,:\,\mbox{problem~\eqref{eq:P} possesses a solution}\right\},\\
 & {\mathcal K}_*:=\left\{\kappa>0\,:\,\mbox{problem~\eqref{eq:P} possesses a minimal solution~$u^\kappa$ such that $u^\kappa\in L^r_{c_*,d_*}$}\right\},\\
 & \kappa^*:=\sup\,{\mathcal K}^*,\quad \kappa_*:=\sup\,{\mathcal K}_*.
\end{split}
\end{equation}
We find $j_*\in\{0,1,2,\dots\}$ such that 
$$
V_j\in L^\infty_{0,-N+2}\quad\mbox{for $j\ge j_*$}
$$
(see Lemma~\ref{Lemma:2.4}), and set $w:=u-V_{j_*}^\kappa$ for $\kappa\in{\mathcal K}^*$. 
Then $u$ is a solution to problem~\eqref{eq:P} if and only if $w$ satisfies an integral equation
\begin{equation}
\label{eq:1.11}
w=\Gamma*\left[\alpha\left((w+U^\kappa_{j_*})^p-(U^\kappa_{j_{*}-1})^{p}\right)\right]\quad\mbox{in}\quad{\mathbb R}^N.
\end{equation}
Applying the same arguments as in \cite{IOS01}, we see that 
\begin{equation}
\label{eq:1.12}
\begin{split}
 & {\mathcal K}_*\subset {\mathcal K}^*=\left\{\kappa>0\,:\,\mbox{problem~\eqref{eq:P} possesses a minimal solution~$u^\kappa$}\right\},\\ 
 & (0,\kappa_*)\subset {\mathcal K}_*,\quad (0,\kappa^*)\subset {\mathcal K}^*,
\end{split}
\end{equation}
if ${\mathcal K}_*\not=\emptyset$ (see Lemma~\ref{Lemma:3.1}). 
Furthermore, 
we apply the contraction mapping theorem in $L^\infty_{0,-N+2}$ to find a function $w\in L^\infty_{0,-N+2}$ satisfying \eqref{eq:1.11} 
for all small enough $\kappa>0$. 
This means that
\begin{equation}
\label{eq:1.13}
{\mathcal K}_*\not=\emptyset
\end{equation}
(see Lemma~\ref{Lemma:3.2}). 
In addition, for any $\kappa\in{\mathcal K}_*$, 
we find the first eigenvalue $\lambda^\kappa$ and its corresponding positive eigenfunction $\varphi^\kappa\in{\mathcal D}^{1,2}$ 
to the linearized eigenvalue problem~\eqref{eq:E} to problem~\eqref{eq:P} at $u^\kappa$ (see Lemma~\ref{Lemma:5.1}). 
Then, thanks to elliptic regularity theorems and the Kelvin transformation, we see that $\varphi^\kappa\in L^\infty_{0,-N+2}$ for $\kappa\in{\mathcal K}_*$ 
(see Lemma~\ref{Lemma:5.2}) and $\lambda^\kappa>1$ for $\kappa\in(0,\kappa_*)$ (see Lemma~\ref{Lemma:5.3}).
We also approximate $\alpha (u^\kappa)^{p-1}$ by functions in $L^\infty_c$ to obtain 
\begin{equation}
\label{eq:1.14}
\int_{{\mathbb R}^N}\alpha p(u^\kappa)^{p-1}\psi^2\,dx
\le\int_{{\mathbb R}^N}|\nabla\psi|^2\,dx,
\quad \psi\in{\mathcal D}^{1,2},
\end{equation}
for all $\kappa\in(0,\kappa^*)$ (see Lemma~\ref{Lemma:5.4}). This leads to that $\kappa^*<\infty$, 
and completes the proof of Theorem~\ref{Theorem:1.1}. 

For the proof of Theorem~\ref{Theorem:1.2}, 
we obtain a uniform energy estimate of $\{(w^\kappa)^\nu\}_{\kappa\in(0,\kappa_*)}$ in $B(0,2)$ for all $\nu\ge 1$ with $\nu^2/(2\nu-1)<p$. 
Then, under a suitable condition on $p$, $\nu$, and $a$, 
we apply elliptic regularity theorems to obtain a uniform $L^\infty(B(0,1))$-estimate of $\{w^\kappa\}_{\kappa\in(0,\kappa_*)}$. 
We also apply the same arguments to the Kelvin transformation of $w^\kappa$. 
Then we obtain a uniform $L^\infty_{0,-N+2}$-estimate of $\{w^\kappa\}_{\kappa\in(0,\kappa_*)}$, and show that $\kappa_*\in{\mathcal K}_*$. 
Furthermore, we prove that $\lambda^{\kappa_*}=1$ and the uniqueness of solutions to problem~\eqref{eq:P} with $\kappa=\kappa_*$. 
Finally, thanks to \eqref{eq:1.14}, we obtain $\kappa_*=\kappa^*$, and complete the proof of Theorem~\ref{Theorem:1.2}. 
\vspace{5pt}

The rest of this paper is organized as follows. 
In Section~2 we obtain preliminary results on the function space $L^q_{\beta,\gamma}$ (see Section~2.1) 
and the functions $\{U_j^\kappa\}$ and $\{V_j^\kappa\}$ (see Section~2.2). 
In Section~3 we prove \eqref{eq:1.12} and \eqref{eq:1.13}.
In Section~4 we study the relation between problem~\eqref{eq:P} and  the Kelvin transformation.  
In Section~5 we study the linearized eigenvalue problem to problem~\eqref{eq:P} at $u^\kappa$, and prove Theorem~\ref{Theorem:1.1}.
Section~6 is devoted to uniform $L^\infty(B(0,1))$-estimates of $\{w^\kappa\}_{\kappa\in(0,\kappa_*)}$ and their Kelvin transformations. 
In Section~7 we complete the proof of Theorem~\ref{Theorem:1.2}. 
Furthermore, we obtain Corollaries~\ref{Corollary:1.1} and \ref{Corollary:1.2}.
\section{Preliminary}
In this section we collect some properties of the function space $L^q_{\beta,\gamma}$. 
Furthermore, we obtain some estimates of the functions $\{U_j^\kappa\}$ and $\{V_j^\kappa\}$. 
In what follows we use $C$ to denote generic positive constants and 
point out that $C$  may take different values  within a calculation.
\subsection{Function space $L^q_{\beta,\gamma}$}
In this subsection we prove the following two lemmas on the function space $L^q_{\beta,\gamma}$. 
\begin{lemma}
\label{Lemma:2.1}
Let $f\in L^q_{\beta,\gamma}$, where $q\in[1,\infty]$ and $\beta$, $\gamma\in{\mathbb R}$. 
\begin{itemize}
  \item[{\rm (1)}] 
  If $q'\le q$, $\beta'\le\beta$, and $\gamma'\ge\gamma$, 
  then there exists $C_1>0$ such that 
  $$
  \|f\|_{L^{q'}_{\beta',\gamma'}}\le C_1\|f\|_{L^q_{\beta,\gamma}}.
  $$
  \item[{\rm (2)}] 
  If $\beta>-N/q$ and $\gamma<-N/q$,
  then there exists $C_2>0$ such that 
  $$
  \|f\|_{L^q}\le C_2\|f\|_{L^q_{\beta,\gamma}}.
  $$
  \item[{\rm (3)}]
  For any $\rho,\sigma\in{\mathbb R}$ and $\tau>0$, set $f_{\rho,\sigma,\tau}(x):=\omega_{\rho,\sigma}(x) |f(x)|^\tau$. 
  Then there exists $C_3>0$ such that 
  $$
  \|f_{\rho,\sigma,\tau}\|_{L^\frac{q}{\tau}_{\tau\beta+\rho,\tau\gamma+\sigma}}\le C_3\|f\|_{L^q_{\beta,\gamma}}^\tau.
  $$
  \item[{\rm (4)}] 
  If $\beta>-N$ and $\gamma\not=-N$, 
  then there exists $C_4>0$ such that 
  $$
  \int_{B(0,R)}|f|\,dx\le C_4\|f\|_{L^q_{\beta,\gamma}}
  \times
  \left\{
  \begin{array}{ll}
  R^{\beta+N} & \mbox{for $R\in(0,1]$},\vspace{5pt}\\
  R^{\max\{\gamma+N,0\}} & \mbox{for $R\in(1,\infty)$}.
  \end{array}
  \right.
  $$
  \item[{\rm (5)}] 
  If $\beta\not=-N$ and $\gamma<-N$, 
  then there exists $C_5>0$ such that 
  $$
  \int_{{\mathbb R}^N\setminus B(0,R)}|f|\,dx\le C_5\|f\|_{L^q_{\beta,\gamma}}
  \times
  \left\{
  \begin{array}{ll}
  R^{\min\{\beta+N,0\}} & \mbox{for $R\in(0,1]$},\vspace{5pt}\\
  R^{\gamma+N} & \mbox{for $R\in(1,\infty)$}.
  \end{array}
  \right.
  $$
\end{itemize}
\end{lemma}
{\bf Proof.}
Assertion~(1) easily follows from H\"older's inequality. 
We prove assertion~(2). 
Assume that $\beta+N/q>0$ and $\gamma+N/q<0$. 
It follows that
\[
\|f\|_{L^q}\le\sum_{i=-\infty}^{\infty}\|f\|_{L^q(A(0,2^i))}
\le\left(\sum_{i=-\infty}^0 2^{\left(\beta+\frac{N}{q}\right)i}+\sum_{i=1}^\infty 2^{\left(\gamma+\frac{N}{q}\right)i}\right)\|f\|_{L^q_{\beta,\gamma}}
\le C\|f\|_{L^q_{\beta,\gamma}}
\]
for $f\in L^q_{\beta,\gamma}$. Thus assertion~(2) holds. 
Assertion~(3) also easily follows from the definition of the norm of $L^q_{\beta,\gamma}$. 

We prove assertion~(4). 
Let $\beta>-N$ and $\gamma\not=-N$. 
Then
$$
\int_{B(0,R)}|f|\,dx=\sum_{i=0}^\infty\int_{A(0,2^{-i}R)}|f|\,dx
\le C\sum_{i=0}^\infty (2^{-i}R)^{\beta+N}\|f\|_{L^q_{\beta,\gamma}}
\le CR^{\beta+N}\|f\|_{L^q_{\beta,\gamma}}
$$
for $R\in(0,1)$. 
Furthermore,
\begin{align*}
\int_{B(0,R)}|f|\,dx & \le\int_{B(0,1/2)}|f|\,dx+\sum_{i=0}^k\int_{A(0,2^{-i}R)}|f|\,dx
\le C\|f\|_{L^q_{\beta,\gamma}}\left(1+\sum_{i=0}^k(2^{-i}R)^{\gamma+N}\right)
\end{align*}
for $R\in[1,\infty)$, where $k$ is the smallest integer satisfying $2^{-(k+1)}R\le 1/2$. 
These imply assertion~(4). 

We prove assertion~(5). 
Let $\beta\not=-N$ and $\gamma<-N$. 
Then
$$
\int_{{\mathbb R}^N\setminus B(0,R)}|f|\,dx=\sum_{i=1}^\infty\int_{A(0,2^i R)}|f|\,dx
\le C\sum_{i=1}^\infty (2^i R)^{\gamma+N}\|f\|_{L^q_{\beta,\gamma}}
\le CR^{\gamma+N}\|f\|_{L^q_{\beta,\gamma}}
$$
for $R\in(1,\infty)$. 
On the other hand, if $\beta>-N$, then
it follows from assertions~(1) and (2) that 
$$
\int_{{\mathbb R}^N\setminus B(0,R)}|f|\,dx\le\|f\|_1\le C\|f\|_{L^1_{\beta,\gamma}}\le C\|f\|_{L^q_{\beta,\gamma}}.
$$
If $\beta<-N$, then 
\begin{align*}
\int_{{\mathbb R}^N\setminus B(0,R)}|f|\,dx & \le \int_{{\mathbb R}^N\setminus B(0,1)}|f|\,dx
+\sum_{i=1}^k\int_{A(0,2^i R)}|f|\,dx\\
 & \le C\left(1+\sum_{i=1}^k(2^i R)^{\beta+N}\right)\|f\|_{L^q_{\beta,\gamma}}
 \le CR^{\beta+N}\|f\|_{L^q_{\beta,\gamma}}
\end{align*}
for $R\in(0,1]$, where $k$ is the smallest integer satisfying $2^kR\ge 1$. 
Then we obtain assertion~(5). 
Thus Lemma~\ref{Lemma:2.1} follows.
$\Box$
\begin{lemma}
\label{Lemma:2.2}
Let $N\ge 3$.
Assume $\beta>-N$, $\beta\not=-2$, $\gamma<-2$, and $\gamma\not=-N$. 
Let $q\in[1,\infty]$ with $q\not=N/2$ and
\begin{equation}
\label{eq:2.1}
\left\{
\begin{array}{ll}
r\in[q,\infty)\quad\mbox{with}\quad
\displaystyle{\frac{1}{q}-\frac{1}{r}<\frac{2}{N}} & \quad\mbox{if}\quad \displaystyle{q<\frac{N}{2}},\vspace{5pt}\\
r=\infty & \quad\mbox{if}\quad \displaystyle{q>\frac{N}{2}}.
\end{array}
\right.
\end{equation}
Then there exists $C>0$ such that 
\begin{equation}
\label{eq:2.2}
\|\Gamma*f\|_{L^r_{\beta_*,\gamma_*}}\le C\|f\|_{L^q_{\beta,\gamma}}
\quad\mbox{for $f\in L^q_{\beta,\gamma}$},
\end{equation}
where $\beta_*:=\min\{\beta+2,0\}$ and $\gamma_*:=\max\{\gamma+2,-N+2\}$.
\end{lemma}
{\bf Proof.}
Let $R>0$. Set 
\begin{align*}
 & g_1(x):=\int_{\Omega_R}\Gamma(x-y)|f(y)|\,dy,\\
 & g_2(x):=\int_{B(0,R/4)}\Gamma(x-y)|f(y)|\,dy,
\quad
g_3(x):=\int_{{\mathbb R}^N\setminus B(0,4R)}\Gamma(x-y)|f(y)|\,dy,
\end{align*}
for $x\in{\mathbb R}^N$, where $\Omega_R:=B(0,4R)\setminus B(0,R/4)$. 
Then $|\Gamma*f|\le g_1+g_2+g_3$ in ${\mathbb R}^N$.

Let $1\le q\le r\le\infty$ and $q<\infty$.  Assume \eqref{eq:2.1}.
Let $s\in[1,\infty)$ be such that 
$$
\frac{1}{s}=1-\frac{1}{q}+\frac{1}{r}=1-\delta\quad\mbox{with}\quad \delta:=\frac{1}{q}-\frac{1}{r}\in\left[0,\frac{2}{N}\right).
$$
Since 
$$
\Gamma(x-y)|f(y)|=\Gamma(x-y)^{\frac{s}{r}}|f(y)|^{\frac{q}{r}}\Gamma(x-y)^{s\left(1-\frac{1}{q}\right)}|f(y)|^{q\delta},
\quad x,y\in{\mathbb R}^N,
$$
by H\"older's inequality we have
\begin{equation}
\label{eq:2.3}
g_1(x)
\le\left(\int_{\Omega_R}\Gamma(x-y)^s|f(y)|^q\,dy\right)^{1/r}
\left(\int_{\Omega_R} \Gamma(x-y)^s\,dy\right)^{1-\frac{1}{q}}
\left(\int_{\Omega_R} |f(y)|^q\,dy\right)^\delta
\end{equation}
for $x\in{\mathbb R}^N$. 
Since $1/s>1-2/N=(N-2)/N$, we have 
$$
\int_{B(0,6R)}\Gamma(y)^s\,dy\le CR^{N-s(N-2)}\quad\mbox{for $R>0$},
$$
which together with \eqref{eq:2.3} implies that
\begin{equation}
\label{eq:2.4}
\begin{split}
 & \|g_1\|_{L^r(A(0,R))}\\
 & \le C\left(\int_{\Omega_R}\left(\int_{A(0,R)}\Gamma(x-y)^s\,dx\right)|f(y)|^q\,dy\right)^{1/r}
 \left(R^{N-s(N-2)}\right)^{1-\frac{1}{q}}
 \|f\|_{L^q(\Omega_R)}^{\delta q}\\
 & \le CR^{-N\left(\frac{1}{q}-\frac{1}{r}\right)+2}\|f\|_{L^q(\Omega_R)}\le CR^{\frac{N}{r}+2}\omega_{\beta,\gamma}(R)\|f\|_{L^q_{\beta,\gamma}}
\end{split}
\end{equation}
for $R>0$ if $r<\infty$. Similarly, we have 
\begin{equation}
\label{eq:2.5}
\|g_1\|_{L^\infty(A(0,R))}\le CR^2\omega_{\beta,\gamma}(R)\|f\|_{L^q_{\beta,\gamma}}\quad\mbox{for $R>0$ if $r=\infty$}.
\end{equation}
By \eqref{eq:2.4} and \eqref{eq:2.5} we obtain 
\begin{equation}
\label{eq:2.6}
R^{-\frac{N}{r}}\|g_1\|_{L^r(A(0,R))}\le CR^2\omega_{\beta,\gamma}(R)\|f\|_{L^q_{\beta,\gamma}}
\quad\mbox{for $R>0$ if $q<\infty$}.
\end{equation}
On the other hand, if $q=\infty$, then
\begin{equation}
\label{eq:2.7}
g_1(x)\le \|f\|_{L^\infty(\Omega_R)}\int_{B(0,6R)}\Gamma(z)\,dz
\le CR^2\|f\|_{L^\infty(\Omega_R)}
\le CR^2\omega_{\beta,\gamma}(R)\|f\|_{L^\infty_{\beta,\gamma}}
\end{equation}
for $x\in A(0,R)$.

Since $\beta>-N$ and $\gamma\not=-N$, 
by Lemma~\ref{Lemma:2.1}~(4)
we have
\begin{equation}
\label{eq:2.8}
\begin{split}
g_2(x)
 & \le C\int_{B(0,R/4)}(|x|-|y|)^{-N+2}|f(y)|\,dy\\
 & \le CR^{-N+2}\int_{B(0,R/4)}|f(y)|\,dy
 \le C\omega_{\beta+2,\gamma_*}(R)\|f\|_{L^q_{\beta,\gamma}}
 \quad\mbox{for $x\in A(0,R)$}.
\end{split}
\end{equation}
Similarly, since 
$\beta-N+2\not=-N$ and $\gamma-N+2<-N$, 
by Lemma~\ref{Lemma:2.1}~(3) and (5)
we obtain
\begin{equation}
\label{eq:2.9}
\begin{split}
 & g_3(x)
 \le C\int_{{\mathbb R}^N\setminus B(0,4R)}(|y|-|x|)^{-N+2}|f(y)|\,dy\\
 & \qquad\quad
 \le C\int_{{\mathbb R}^N\setminus B(0,4R)}|y|^{-N+2}|f(y)|\,dy
 \le C\omega_{\beta_*,\gamma+2}(R)\|\tilde{f}\|_{L^q_{\beta-N+2,\gamma-N+2}}\\
 & \qquad\quad
 \le C\omega_{\beta_*,\gamma+2}(R)\|f\|_{L^q_{\beta,\gamma}}
  \quad\mbox{for $x\in A(0,R)$},
\end{split}
\end{equation}
where $\tilde{f}(x):=|x|^{-N+2}|f(x)|$.  
Therefore we deduce from \eqref{eq:2.6}, \eqref{eq:2.7}, \eqref{eq:2.8}, and \eqref{eq:2.9} that 
$$
R^{-\frac{N}{r}}\|\Gamma*f\|_{L^r(A(0,R))} \le C\omega_{\beta_*,\gamma_*}(R)\|f\|_{L^q_{\beta,\gamma}}
\quad\mbox{for $R>0$}.
$$
This implies \eqref{eq:2.2}, and Lemma~\ref{Lemma:2.2} follows. 
$\Box$
\subsection{Approximate solutions}
Assume \eqref{eq:1.3}. 
Let $\mu\in{\mathcal M}_+$ be such that $\Gamma*\mu\in L^r_{c,d}$ for some $r\in(p,\infty]$. 
We obtain some estimates of the functions $\{U_j^\kappa\}_{j=0}^\infty$ and $\{V_j^\kappa\}_{j=0}^\infty$.
\begin{lemma}
\label{Lemma:2.3}
Assume the same conditions as in Theorem~{\rm \ref{Theorem:1.1}}. Let $0<\kappa\le\kappa'$ and $j=0,1,\ldots$. 
Then
\begin{align}
\label{eq:2.10}
 & \frac{\kappa'}{\kappa}U_j^\kappa(x)\le U_j^{\kappa'}(x)\le\left(\frac{\kappa'}{\kappa}\right)^{p^j}U^\kappa_j(x),\\
\label{eq:2.11}
 & \left(\frac{\kappa'}{\kappa}\right)^{j(p-1)+1}V_j^\kappa(x)\le V_j^{\kappa'}(x)\le\left(\frac{\kappa'}{\kappa}\right)^{p^j}V_j^\kappa(x),
\end{align}
for a.a.~$x\in{\mathbb R}^N$. 
Furthermore, there exists $C>0$ such that
\begin{equation}
\label{eq:2.12}
\begin{split}
 & 0<\varepsilon U_j^\kappa(x)\le U_j^{(1+\varepsilon)\kappa}(x)-U_j^\kappa(x) \le C\varepsilon U_j^\kappa(x),\\
 & 0<\varepsilon V_j^\kappa(x) \le V_j^{(1+\varepsilon)\kappa}(x)-V_j^\kappa(x)\le C\varepsilon V_j^\kappa(x),
\end{split}
\end{equation}
for a.a.~$x\in{\mathbb R}^N$ and $0<\varepsilon<1$.
\end{lemma}
{\bf Proof.}
Let $0<\kappa\le\kappa'$. 
By induction we prove \eqref{eq:2.10}. 
By the definition of $U_0^\kappa$ we easily obtain \eqref{eq:2.10} for $j=0$. 
Assume that \eqref{eq:2.10} holds for some $j=j_0\in\{0,1,2,\dots\}$. 
Then 
\begin{align*}
U_{j_0+1}^{\kappa'}
 & \ge\Gamma*\alpha\left(\frac{\kappa'}{\kappa}U_{j_0}^\kappa\right)^p+\kappa' (\Gamma*\mu)
\ge\frac{\kappa'}{\kappa}
\left[\Gamma*\alpha (U_{j_0}^\kappa)^p+\kappa(\Gamma*\mu)\right]\\
 & =\frac{\kappa'}{\kappa} U_{j_0+1}^\kappa,\\
U_{j_0+1}^{\kappa'}
 & \le\Gamma*\alpha\biggr(\left(\frac{\kappa'}{\kappa}\right)^{p^{j_0}}U^\kappa_{j_0}\biggr)^p+\kappa' (\Gamma*\mu)
\le\left(\frac{\kappa'}{\kappa}\right)^{p^{j_0+1}}
\left[\Gamma*\alpha (U_{j_0}^\kappa)^p+\kappa(\Gamma*\mu)\right]\\
 & =\left(\frac{\kappa'}{\kappa}\right)^{p^{j_0+1}} U_{j_0+1}^\kappa.
\end{align*}
These imply that \eqref{eq:2.10} holds for $j=j_0+1$. Thus \eqref{eq:2.10} holds for all $j=0,1,2,\dots$.

Next, we prove \eqref{eq:2.11}. 
Since 
$$
V_0^{\kappa'}=\kappa' [\Gamma*\mu](x)=\frac{\kappa'}{\kappa}V_0^\kappa,
\quad
V_1^{\kappa'}=\Gamma*(\alpha U_0^{\kappa'})^p
=\left(\frac{\kappa'}{\kappa}\right)^pV_1^\kappa,
$$
we have \eqref{eq:2.11} for $j=0,1$. 
Assume that \eqref{eq:2.11} holds for some $j=j_0\in\{1,2,\dots\}$. 
Then, for any $t\in(0,1)$, by \eqref{eq:2.10} we have
\begin{align*}
p((1-t)U_{j_0-1}^{\kappa'}+tU_{j_0}^{\kappa'})^{p-1}V_{j_0}^{\kappa'} & 
\ge\left(\frac{\kappa'}{\kappa}\right)^{p-1}p((1-t)U_{j_0-1}^\kappa+tU_{j_0}^\kappa)^{p-1}
\left(\frac{\kappa'}{\kappa}\right)^{(p-1)j_0+1}V_{j_0}^\kappa\\
&=\left(\frac{\kappa'}{\kappa}\right)^{(p-1)(j_0+1)+1}p((1-t)U_{j_0-1}^\kappa+tU_{j_0}^\kappa)^{p-1}V_{j_0}^\kappa,\\
p((1-t)U_{j_0-1}^{\kappa'}+tU_{j_0}^{\kappa'})^{p-1}V_{j_0}^{\kappa'} & 
\le\left(\frac{\kappa'}{\kappa}\right)^{(p-1)p^{j_0}}p((1-t)U^\kappa_{j_0-1}+tU^\kappa_{j_0})^{p-1}
\left(\frac{\kappa'}{\kappa}\right)^{p^{j_0}}V^\kappa_{j_0}\\
&=\left(\frac{\kappa'}{\kappa}\right)^{p^{j_0+1}}p((1-t)U^\kappa_{j_0-1}+tU^\kappa_{j_0})^{p-1}V^\kappa_{j_0}.
\end{align*}
Since
\[
V^\kappa_{j_0+1}=\Gamma*\left[\alpha((U^\kappa_{j_0})^p-(U^\kappa_{j_0-1})^p)\right]
=\Gamma*\left[\alpha\int_0^1p((1-t)U^\kappa_{j_0-1}+tU^\kappa_{j_0})^{p-1}V^\kappa_{j_0}\,dt\right],
\]
we obtain
\[
\left(\frac{\kappa'}{\kappa}\right)^{(j_0+1)(p-1)+1}V^\kappa_{j_0+1}
\le V^{\kappa'}_{j_0+1}\le\left(\frac{\kappa'}{\kappa}\right)^{p^{j_0+1}}V^\kappa_{j_0+1}.
\]
This implies that \eqref{eq:2.11} holds for $j=j_0+1$. Thus \eqref{eq:2.11} holds for all $j=0,1,2,\dots$.
Furthermore, relation~\eqref{eq:2.12} follows from \eqref{eq:2.10} and \eqref{eq:2.11}. 
Thus Lemma~\ref{Lemma:2.3} follows.
$\Box$\vspace{5pt}

Let $r\in(p,\infty)$ and assume \eqref{eq:1.3}. Then
\begin{equation}\label{eq:2.13}
q:=\frac{r}{p-1}>\frac{N}{2},\quad \overline{c}:=(p-1)c_*+a>-2,\quad \overline{d}:=(p-1)d_*+b<-2,
\end{equation}
where $c_*$ and $d_*$ are as in Theorem~\ref{Theorem:1.1}. 
We find $r_*\in(1,\infty)$ such that 
\begin{equation}
\label{eq:2.14}
\max\left\{\frac{N}{2},\frac{r}{r-1}\right\}<r_*<q,
\quad
\frac{1}{r}\not\in\left\{j\left(\frac{2}{N}-\frac{1}{r_*}\right)\,:\,j=0,1,2,\dots\right\}.
\end{equation}
Define a sequence $\{r_j\}_{j=0}^\infty\subset(p,\infty]$ by 
$$
\frac{1}{r_j}:=\max\left\{\frac{1}{r}-j\left(\frac{2}{N}-\frac{1}{r_*}\right),0\right\}.
$$
Also, by \eqref{eq:2.13} we find $\hat{c}$ and $\hat{d}$ such that
\begin{equation}
\label{eq:2.15}
\overline{c}>\hat{c}>-2,\quad
\overline{d}<\hat{d}<-2.
\end{equation}
Define sequences $\{c_j\}_{j=0}^\infty$, $\{d_j\}_{j=0}^\infty$ by
\begin{equation}
\label{eq:2.16}
c_j:=\min\{c+j(\hat{c}+2),0\},
\qquad
d_j:=\max\{d+j(\hat{d}+2),-N+2\}.
\end{equation}
By \eqref{eq:2.14} and \eqref{eq:2.15} there exists $j_*\in\{1,2,\dots\}$ such that 
\begin{equation}
\label{eq:2.17}
r_j=\infty,\qquad
c_j=0,\qquad
d_j=-2+N,\qquad
\mbox{for all $j\ge j_*$}.
\end{equation}
\begin{Lemma}
\label{Lemma:2.4}
Assume the same conditions as in Theorem~{\rm\ref{Theorem:1.1}}. 
Let $K>0$. 
For any $j=0,1,\ldots$, 
there exists $C>0$ such that
\begin{equation}
\label{eq:2.18}
\|U^\kappa_j\|_{L^r_{c_*,d_*}}\le C\kappa,
\qquad
\|V^\kappa_j\|_{L^{r_j}_{c_j,d_j}}\le C\kappa^{(p-1)j+1},
\end{equation}
for $\kappa\in(0,K)$.
\end{Lemma}
{\bf Proof.} 
Since $U_0^\kappa=V_0^\kappa=\kappa\Gamma*\mu$ in ${\mathbb R}^N$, 
by Lemma~\ref{Lemma:2.1}~(1) and \eqref{eq:1.4} we have \eqref{eq:2.18} with $j=0$. 
Assume that \eqref{eq:2.18} holds for some $j=j_0\in\{0,1,2,\dots\}$. 
Notice that \eqref{eq:1.1} implies $\alpha(x)\le C\omega_{a,b}(x)$ for some constant $C>0$. 
By Lemma~\ref{Lemma:2.1}~(3) we obtain
$$
\left\|\alpha(U^\kappa_{j_0})^{p-1}\right\|_{L^q_{\overline{c},\overline{d}}}\le\|U^\kappa_{j_0}\|_{L^r_{c_*,d_*}}^{p-1}\le C\kappa^{p-1}.
$$
Then we observe from H\"older inequality that
\begin{equation}
\label{eq:2.19}
\begin{split}
\left\|\alpha(U^\kappa_{j_0})^p\right\|_{L^{q'}_{c+\overline{c},d+\overline{d}}}
 & \le\left\|\alpha(U^\kappa_{j_0})^{p-1}\right\|_{L^q_{\overline{c},\overline{d}}}\|U^\kappa_{j_0}\|_{L^r_{c_*,d_*}}
\le C\kappa^p,\\
\left\|\alpha(U^\kappa_{j_0})^{p-1}V^\kappa_{j_0}\right\|_{L^{r'}_{c_{j_0}+\overline{c},d_{j_0}+\overline{d}}}
 & \le\left\|\alpha(U^\kappa_{j_0})^{p-1}\right\|_{L^q_{\overline{c},\overline{d}}}\|V^\kappa_{j_0}\|_{L^{r_{j_0}}_{c_{j_0},d_{j_0}}}\le C\kappa^{(p-1)(j_0+1)+1},
\end{split}
\end{equation}
where
$$
\frac{1}{q'}=\frac{1}{q}+\frac{1}{r},
\qquad
\frac{1}{r'}=\frac{1}{q}+\frac{1}{r_{j_0}}.
$$
It follows from \eqref{eq:2.13}, \eqref{eq:2.14}, \eqref{eq:2.15}, and \eqref{eq:2.16} that
\begin{equation*}
\begin{aligned}
\frac{1}{q'}-\frac{1}{r}=\frac{1}{q}<\frac{2}{N}&,
\quad
\frac{1}{r'}-\frac{1}{r_{j_0+1}}\le\frac{1}{q}-\frac{1}{r_*}+\frac{2}{N}<\frac{2}{N},\\
(c+\overline{c})-c=\overline{c}>-2&,\quad (c_{j_0}+\overline{c})-c_{j_0+1}\ge \overline{c}-\hat{c}-2>-2,\\
(d+\overline{d})-d=\overline{d}<-2&,\quad (d_{j_0}+\hat{d})-d_{j_0+1}\le\overline{d}-\hat{d}-2<-2.
\end{aligned}
\end{equation*}
Let $K>0$. 
Then, by Lemma~\ref{Lemma:2.2}, \eqref{eq:1.8}, \eqref{eq:1.9}, and \eqref{eq:2.19} we have 
\begin{align*}
\|U_{j_0+1}^\kappa\|_{L^r_{c_*,d_*}}
 & \le \left\|\Gamma*\left[\alpha (U_{j_0}^\kappa)^p\right]\right\|_{L^r_{c_*,d_*}}+\kappa\|\Gamma*\mu\|_{L^r_{c,d}}
\le C\kappa,\\ 
\|V^\kappa_{j_0+1}\|_{L^{r_{j_0+1}}_{c_{j_0+1},d_{j_0+1}}}
 & \le p\left\|\Gamma*\left[\alpha(U^\kappa_{j_0})^{p-1}V^\kappa_{j_0}\right]\right\|_{L^{r_{j_0+1}}_{c_{j_0+1},d_{j_0+1}}}\le C\kappa^{(p-1)(j_0+1)+1},
\end{align*}
for $\kappa\in(0,K)$. Thus \eqref{eq:2.18} holds for some $j=j_0+1$. 
Therefore we obtain \eqref{eq:2.18} for $j=0,1,2,\dots$, and the proof is complete. 
$\Box$
\section{Solutions to problem~\eqref{eq:P}}
Assume the same conditions as in Theorem~\ref{Theorem:1.1}. 
We denote by $u^\kappa$ the (unique) minimal solution to problem~\eqref{eq:P}. 
We first show that, for problem~\eqref{eq:P}, 
the existence of supersolutions implies the existence of the minimal solution~$u^\kappa$.
\begin{lemma}
\label{Lemma:3.1}
Let $\kappa>0$, and 
assume that there exists a supersolution~$v$ to problem~\eqref{eq:P}. 
Then there exists a minimal solution $u^\kappa$ to problem~\eqref{eq:P} and 
$$
u^\kappa(x)=\lim_{j\to\infty}U_j^\kappa(x)\le v(x)
$$ 
for a.a.~$x\in{\mathbb R}^N$. 
Furthermore, there exists $C>0$, independent of $\kappa>0$, such that 
\begin{equation}
\label{eq:3.1}
u^\kappa(x)\ge C\kappa(1+|x|)^{-N+2}
\end{equation}
for a.a.~$x\in{\mathbb R}^N$.
\end{lemma}
{\bf Proof.}
Let $v$ be a supersolution to problem~\eqref{eq:P}. 
Due to \eqref{eq:1.8}, 
by induction, for any $j\ge 0$, we see that 
$U^\kappa_j(x)\le v(x)$ for a.a.~$x\in{\mathbb R}^N$. 
This together with \eqref{eq:1.9} implies that, 
for a.a.~$x\in{\mathbb R}^N$, 
the limit $U(x):=\lim_{j\to\infty}U^\kappa_j(x)$ exists and $U(x)\le v(x)$. 
Then we observe from \eqref{eq:1.8} that $U$ is a solution to problem~\eqref{eq:P}. 
Furthermore, we see that $U$ is a minimal solution to problem~\eqref{eq:P}, 
that is, $U=u^\kappa$. 
In addition, it follows that
$$
u^\kappa(x)\ge U_0^\kappa(x)=\kappa\int_{{\mathbb R}^N}\Gamma(x-y)\,d\mu(y)
$$
for a.a.~$x\in{\mathbb R}^N$. 
Since $\mu\in{\mathcal M}_+$, there exists a bounded measurable set $D$ such that $\mu(D)>0$. 
Then 
$$
u^\kappa(x)\ge \kappa\mu(D)\min_{y\in D}\Gamma(x-y)\ge C\kappa(1+|x|)^{-N+2}
$$
for a.a.~$x\in{\mathbb R}^N$. Thus \eqref{eq:3.1} holds, and the proof is complete. 
$\Box$\vspace{5pt}

Let ${\mathcal K}^*$, ${\mathcal K}_*$, $\kappa^*$, and $\kappa_*$ be as in \eqref{eq:1.10}.
Then relation~\eqref{eq:1.12} follows from Lemma~\ref{Lemma:3.1}. 
Furthermore, we have:
\begin{Lemma}
\label{Lemma:3.2}
Assume the same conditions as in Theorem~{\rm\ref{Theorem:1.1}}. 
Then $\kappa_*>0$. 
\end{Lemma}
{\bf Proof.}
For any $v\in L^\infty_{0,-N+2}$, we define
\[
\Phi[v]:=\Gamma*\left[\alpha\left(\left(v_++U^\kappa_{j_*}\right)^p-\left(U^\kappa_{j_{*}-1}\right)^p\right)\right],
\]
where $v_+:=\max\{v,0\}$. 
Since $q>N/2$, $\overline{c}>-2$, and $\overline{d}<-2$ (see \eqref{eq:2.13}), 
by \eqref{eq:1.1} we apply Lemmas~\ref{Lemma:2.1}, \ref{Lemma:2.2}, and \ref{Lemma:2.4} to obtain 
\begin{equation}
\label{eq:3.2}
\begin{aligned}
\|\Phi[v]\|_{L^\infty_{0,-N+2}}&\le p\left\|\Gamma*\left[\alpha \left(v_++U^\kappa_{j_*}\right)^{p-1}V^\kappa_{j_*}\right]\right\|_{L^\infty_{0,-N+2}}\\
&\le C\left\|\alpha\left(v_++U^\kappa_{j_*}\right)^{p-1}V^\kappa_{j_*}\right\|_{L^{q}_{\overline{c},\overline{d}-N+2}}
\le C\left\|v_++U^\kappa_{j_*}\right\|_{L^r_{c_*,d_*}}^{p-1}\|V^\kappa_{j_*}\|_{L^\infty_{0,-N+2}}\\
&\le C\kappa^{(p-1)j_*+1}(\|v\|_{L^\infty_{0,-N+2}}+\kappa)^{p-1}
\end{aligned}
\end{equation}
for all $v\in L^\infty_{0,-N+2}$ 
and
\begin{equation}
\begin{aligned}
\|\Phi[v_1]-\Phi[v_2]\|_{L^\infty_{0,-N+2}} & 
\le p\left\|\Gamma*\left[\alpha(\max\{(v_1)_+,(v_2)_+\}+U^\kappa_{j_*})^{p-1}(v_1-v_2)\right]\right\|_{L^\infty_{0,-N+2}}\\
&\le C\|\alpha(\max\{(v_1)_+,(v_2)_+\}+U^\kappa_{j_*})^{p-1}(v_1-v_2)\|_{L^{q}_{\overline{c},\overline{d}-N+2}}\\
&\le C\|\max\{(v_1)_+,(v_2)_+\}+U^\kappa_{j_*}\|_{L^r_{c_*,d_*}}^{p-1}\|v_1-v_2\|_{L^\infty_{0,-N+2}}\\
&\le C\left(\max\{\|v_1\|_{L^\infty_{0,-N+2}},\|v_2\|_{L^\infty_{0,-N+2}}\}+\kappa\right)^{p-1}\|v_1-v_2\|_{L^\infty_{0,-N+2}}
\end{aligned}
\label{eq:3.3}\end{equation}
for all $v_1$, $v_2\in L^\infty_{0,-N+2}$. 
By \eqref{eq:3.2} and \eqref{eq:3.3}, 
for any small enough $\kappa>0$, 
we see that $\Phi$ is a contraction map on the ball 
$$
{\mathcal B}_\kappa:=\{v\in L^\infty_{0,-N+2}\,:\, \|v\|_{L^\infty_{0,-N+2}}\le\kappa\}. 
$$
Therefore, by the contraction mapping theorem we find $w\in{\mathcal B}_\kappa$ such that $w=\Phi[w]$. 
Then $w$ satisfies 
$$
w>0\quad\mbox{in}\quad{\mathbb R}^N,\qquad
w=\Gamma*\left[\alpha((w+U^\kappa_{j_*})^p-(U^\kappa_{j_{*}-1})^{p})\right]\quad\mbox{in}\quad{\mathbb R}^N.
$$
Therefore, setting $u=w+U_{j_*}^\kappa$, we see that $u$ is a solution to problem~\eqref{eq:P}. 
Furthermore, we observe from Lemma~\ref{Lemma:2.4} that $u\in L^r_{c_*,d_*}$. 
Thus Lemma~\ref{Lemma:3.2} follows.
$\Box$\vspace{5pt}

By \eqref{eq:1.12}, for any $\kappa\in{\mathcal K}^*$, 
there exists a minimal solution $u^\kappa$ to problem~\eqref{eq:P}. 
Set 
\begin{equation}
\label{eq:3.4}
w^\kappa(x):=u^\kappa(x)-U_{j_*}^\kappa(x)\quad\mbox{for a.a.~$x\in{\mathbb R}^N$}.
\end{equation}
Then $w^\kappa$ satisfies 
\begin{equation}
\label{eq:3.5} 
w^\kappa(x)=\Gamma*\left[\alpha((w^\kappa+U^\kappa_{j_*})^p-(U^\kappa_{j_*-1})^p)\right]
\le \Gamma*\left[p\alpha(u^\kappa)^{p-1}(w^\kappa+V^\kappa_{j_*})\right]
\end{equation}
for a.a.~$x\in{\mathbb R}^N$. 
Furthermore, 
by Lemmas~\ref{Lemma:2.3} and \ref{Lemma:3.1} we see that, 
if $0<\kappa<\kappa'$ and $\kappa'\in{\mathcal K}^*$, then 
\begin{equation}
\label{eq:3.6}
\begin{split}
w^{\kappa'}(x) & =u^{\kappa'}(x)-U_{j_*}^{\kappa'}(x)
=\lim_{j\to\infty}\left(U_j^{\kappa'}(x)-U_{j_*}^{\kappa'}(x)\right)\\
 & =\lim_{j\to\infty} \left(V_j^{\kappa'}(x)+V_{j-1}^{\kappa'}(x)+\cdots+V_{j_*+1}^{\kappa'}(x)\right)\\
 & >\lim_{j\to\infty} \left(V_j^\kappa(x)+V_{j-1}^\kappa(x)+\cdots+V_{j_*+1}^\kappa(x)\right)\\
 & =\lim_{j\to\infty}\left(U_j^\kappa(x)-U_{j_*}^\kappa(x)\right)
 =u^\kappa(x)-U_{j_*}^\kappa(x)
=w^\kappa(x)
\end{split}
\end{equation}
for a.a.~$x\in{\mathbb R}^N$. 
In addition, we have:
\begin{Lemma}
\label{Lemma:3.3}
Assume the same conditions as in Theorem~{\rm\ref{Theorem:1.1}}. 
Let $\kappa\in{\mathcal K}_*$ and let $w^\kappa$ be as in~\eqref{eq:3.4}. 
Then 
$w^\kappa\in {\mathcal D}^{1,2}\cap L^\infty_{0,-N+2}$ and 
\begin{equation}
\label{eq:3.7}
\int_{{\mathbb R}^N}\nabla w^\kappa\cdot\nabla\psi\,dx=\int_{{\mathbb R}^N}\alpha((w^\kappa+U^\kappa_{j_*})^p-(U^\kappa_{j_*-1})^p)\psi\,dx
\end{equation}
holds for all $\psi\in {\mathcal D}^{1,2}$.
\end{Lemma}
{\bf Proof.} 
Let $\kappa\in{\mathcal K}_*$. 
It follows from $0\le w^\kappa\le u^\kappa$ 
that $w^\kappa\in L^r_{c_*,d_*}$. 
By Lemma~\ref{Lemma:2.4} and \eqref{eq:2.17} 
we see that $V_{j_*}^\kappa\in L^{r_{j_*}}_{c_{j_*},d_{j_*}}=L^\infty_{0,-N+2}$. 
Then, applying the same argument as in the proof of \eqref{eq:2.18} to \eqref{eq:3.5}, 
we see that
\begin{equation}
\label{eq:3.8}
\alpha(u^\kappa)^{p-1}\in L^{q}_{\overline{c},\overline{d}},\qquad
w^\kappa\in L^{r_j}_{r_j,d_j}\quad\mbox{for $j=1,2,\dots$}.
\end{equation}
We deduce that 
$$
w^\kappa\in L^\infty_{0,-N+2},\quad 
\alpha((w^\kappa+U^\kappa_{j_*})^p-(U^\kappa_{j_*-1})^p)
\le p\alpha(u^\kappa)^{p-1}(w^\kappa+V^\kappa_{j_*})\in L^{q}_{\overline{c},\overline{d}-N+2}.
$$
On the other hand, 
it follows from \eqref{eq:2.13}
that $q>N/2$, $\overline{c}>-2$, and $\overline{d}-N+2<-N$.
Then, by Lemma~\ref{Lemma:2.1}~(2) we have 
$$
\alpha((w^\kappa+U^\kappa_{j_*})^p-(U^\kappa_{j_*-1})^p)\in L^s({\mathbb R}^N)\quad\mbox{with $s\in\displaystyle{\left[1,\frac{N}{2}\right]}$}.
$$
Therefore we deduce from \cite{GT}*{Section~9.4} that
$$
\nabla^2 w^\kappa\in L^s({\mathbb R}^N)
\quad\mbox{with $s\in\displaystyle{\left(1,\frac{N}{2}\right]}$},
\qquad
w^\kappa\in W^{2,N/2}_{{\rm loc}}({\mathbb R}^N),
$$ 
and $w^\kappa$ satisfies \eqref{eq:3.7}.
In particular, 
the Sobolev imbedding theorem implies that $\nabla w^\kappa\in L^2({\mathbb R}^N)$.
Thus Lemma~\ref{Lemma:3.3} follows. 
$\Box$
\section{Dual problem via the Kelvin transform}
Let $T$ be a map from ${\mathbb R}^N\setminus\{0\}$ into ${\mathbb R}^N$ such that 
$$
Tx:=\frac{x}{|x|^2}\quad\mbox{for $x\in{\mathbb R}^N\setminus\{0\}$}.
$$
For any measurable function $f$ in ${\mathbb R}^N$, we define the Kelvin transformation $f^\sharp$ of $f$ by 
$$
f^\sharp(x):=|x|^{-N+2}f(Tx)\quad\mbox{for a.a.~$x\in{\mathbb R}^N$}.
$$
\begin{lemma}
\label{Lemma:4.1}
Assume the same conditions as in Theorem~{\rm \ref{Theorem:1.1}}. 
Let $\kappa>0$, and let $u$ be a solution to problem~\eqref{eq:P}. 
Then the Kelvin transformation $u^\sharp$ of $u$ is a solution to the integral equation 
\begin{equation}
\tag{I}
\label{eq:I}
v(x)=\int_{{\mathbb R}^N}\Gamma(x-y)\beta(y)y(y)^p\,dy+\kappa(\Gamma*\mu)^\sharp(x)\quad\mbox{in}\quad{\mathbb R}^N,
\end{equation}
where $\beta(x)=|x|^{(N-2)(p-1)-4}\alpha(Tx)$. 
Here 
$$
\limsup_{|x|\to 0}|x|^{-a^\sharp}\beta(x)<\infty,\quad \limsup_{|x|\to\infty}|x|^{-b^\sharp}\beta(x)<\infty,
\quad (\Gamma*\mu)^\sharp\in L^r_{c^\sharp,d^\sharp},
$$
where 
\begin{equation*}
\begin{array}{ll}
a^\sharp:=(N-2)(p-1)-4-b>-2,\quad 
 & b^\sharp:=(N-2)(p-1)-4-a,\vspace{5pt}\\
c^\sharp:=-N+2-d,\quad & d^\sharp:=-N+2-c.
\end{array}
\end{equation*} 
Furthermore, 
$$
p>\frac{N+b^\sharp}{N-2},\qquad c^\sharp>-\frac{a^\sharp+2}{p-1},\qquad 
d^\sharp<-\frac{b^\sharp+2}{p-1}.
$$
\end{lemma}
{\bf Proof.}
Since 
$$
|Tx-y|^2=|Tx|^2|x-Ty|^2|y|^2,\quad x,y\in{\mathbb R}^N\setminus\{0\},
$$
it follows that 
\begin{equation*}
\begin{split}
|x|^{-N+2}\Gamma(Tx-y) & =\frac{1}{N(N-2)\omega_N}|x|^{-N+2}|Tx-y|^{-N+2}\\
 & =\frac{1}{N(N-2)\omega_N}|x-Ty|^{-N+2}|y|^{-N+2}
=\Gamma(x-Ty)|Ty|^{N-2}
\end{split}
\end{equation*}
for $x$, $y\in{\mathbb R}^N\setminus\{0\}$ with $x\not=y$. 
Then, by Definition~\ref{Definition:1.1}~(1) we have
\begin{align*}
u^\sharp(x) 
 & =\int_{{\mathbb R}^N}\Gamma(x-Ty)|Ty|^{N-2}\alpha(y)u(y)^p\,dy+\kappa(\Gamma*\mu)^\sharp(x)\\
 & =\int_{{\mathbb R}^N}\Gamma(x-y)|y|^{-N-2}\alpha(Ty)u(Ty)^p\,dy+\kappa(\Gamma*\mu)^\sharp(x)\\
 & =\int_{{\mathbb R}^N}\Gamma(x-y)\beta(y)u^\sharp(y)^p\,dy+\kappa(\Gamma*\mu)^\sharp(x)
\end{align*}
for a.a.~$x\in{\mathbb R}^N$. 
This means that $u^\sharp$ is a solution to integral equation~\eqref{eq:I}. 
Furthermore,  
$(\Gamma*\mu)^\sharp\in L^r_{c^\sharp,d^\sharp}$
and 
\begin{equation*}
\begin{split}
a^\sharp & =(N-2)(p-1)-4-b>(N-2)\frac{2+b}{N-2}-4-b=-2,\\
\frac{N+b^\sharp}{N-2}
 & =\frac{Np-2p-2-a}{N-2}
=p-\frac{a+2}{N-2}<p,\\
-\frac{a^\sharp+2}{p-1} & =-N+2+\frac{b+2}{p-1}
<-N+2-d=c^\sharp,\\
-\frac{b^\sharp+2}{p-1} & =-N+2+\frac{a+2}{p-1}
>-N+2-c=d^\sharp.
\end{split}
\end{equation*}
Thus Lemma~\ref{Lemma:4.1} follows.
$\Box$\vspace{5pt}

We easily observe from Lemma~\ref{Lemma:4.1} that $(u^\kappa)^\sharp$ is a minimal solution to integral equation~\eqref{eq:I}. 
Furthermore, we see that
\begin{equation}
\label{eq:4.1}
\begin{array}{ll}
(U_j^\kappa)^\sharp(x)=\left[\Gamma*\beta\left((U_{j-1}^\kappa)^\sharp\right)^p\right](x)
+\kappa (\Gamma*\mu)^\sharp(x),\qquad  & j=0,1,2,\dots,\vspace{5pt}\\
(V_j^\kappa)^\sharp(x)=(U_j^\kappa)^\sharp(x)-(U_{j-1}^\kappa)^\sharp(x),\qquad & j=0,1,2,\dots,\vspace{5pt}\\
(w^\kappa)^\sharp(x)=(u^\kappa-U_{j_*}^\kappa)^\sharp(x)
=(u^\kappa)^\sharp(x)-(U_{j_*}^\kappa)^\sharp(x), & 
\end{array}
\end{equation}
for a.a.~$x\in{\mathbb R}^N$ and 
$(V_j^\kappa)^\sharp\in L^\infty_{0,-N+2}$ for $j\ge j_*$.
\vspace{5pt}
%

At the end of this section we state a lemma on the relation between weak solutions and the Kelvin transformation.
\begin{lemma}
\label{Lemma:4.2}
{\rm (1)} 
Let $\varphi$, $\psi\in{\mathcal D}^{1,2}$. Then $\varphi^\sharp$, $\psi^\sharp\in {\mathcal D}^{1,2}$ and 
$$
\int_{{\mathbb R}^N}\nabla\varphi^\sharp\cdot\nabla\psi^\sharp\,dx=\int_{{\mathbb R}^N}\nabla\varphi\cdot\nabla\psi\,dx.
$$
{\rm (2)}
Let $f$ be a measurable function in ${\mathbb R}^N$.
Assume that $\varphi\in {\mathcal D}^{1,2}$ satisfies
\begin{equation}
\label{eq:4.2}
\int_{{\mathbb R}^N}\nabla\varphi\cdot\nabla\psi\,dx=\int_{{\mathbb R}^N}f\psi\,dx
\quad \mbox{for all $\psi\in{\mathcal D}^{1,2}$}.
\end{equation}
Then 
$$
\int_{{\mathbb R}^N}\nabla\varphi^\sharp\cdot\nabla\psi\,dx=\int_{{\mathbb R}^N}|x|^{-4}f^\sharp \psi\,dx
\quad \mbox{for all $\psi\in{\mathcal D}^{1,2}$}.
$$
\end{lemma}
{\bf Proof.}
It follows that 
\begin{align*}
\int_{{\mathbb R}^N}\nabla\varphi^\sharp\cdot\nabla\psi^\sharp\,dx & =\int_{{\mathbb R}^N}(-\Delta\varphi^\sharp)\psi^\sharp\,dx
 =\int_{{\mathbb R}^N}|x|^{-2N}(-\Delta\varphi)(Tx)\psi(Tx)\,dx\\
 & =\int_{{\mathbb R}^N}(-\Delta\varphi)\psi\,dx=\int_{{\mathbb R}^N}\nabla\varphi\cdot\nabla\psi\,dx
\end{align*}
for $\varphi$, $\psi\in C_c^\infty({\mathbb R}^N)$. 
This implies assertion~(1). 
Furthermore, we observe from assertion~(1) and \eqref{eq:4.2} that
\begin{align*}
\int_{{\mathbb R}^N}\nabla\varphi^\sharp \cdot\nabla\psi\,dx
 & =\int_{{\mathbb R}^N}\nabla\varphi\cdot\nabla\psi^\sharp\,dx
=\int_{{\mathbb R}^N}f\psi^\sharp\,dx\\
 & =\int_{{\mathbb R}^N}|x|^{-2N}f(Tx)\psi^\sharp(Tx)\,dx=\int_{{\mathbb R}^N}|x|^{-4}f^\sharp \psi\,dx
\end{align*}
for all $\psi\in{\mathcal D}^{1,2}$.
This implies assertion~(2), and Lemma~\ref{Lemma:4.2} follows. 
$\Box$
\section{Linearized eigenvalue problems}
Let $\kappa\in{\mathcal K}_*$. 
It follows from \eqref{eq:3.8} that $\alpha (u^\kappa)^{p-1}\in L^q_{\overline{c},\overline{d}}$, 
which together with Lemma~\ref{Lemma:2.1}~(2) and \eqref{eq:2.13} implies that 
\begin{equation}
\label{eq:5.1}
\alpha (u^\kappa)^{p-1}\in L^{q'}({\mathbb R}^N)\cap L^{\frac{N}{2}}({\mathbb R}^N)
\quad\mbox{for some $q'>N/2$}. 
\end{equation}
Then, applying the same argument as in \cite{NS01}*{Lemma~B.2}, 
we have the following lemma on the linearized eigenvalue problem to problem~\eqref{eq:P} at $u=u^\kappa$. 
\begin{lemma}
\label{Lemma:5.1}
Assume the same conditions as in Theorem~{\rm\ref{Theorem:1.1}}. 
Let $\kappa\in{\mathcal K}_*$. 
Consider the eigenvalue problem 
\begin{equation}
\label{eq:E}
\tag{$\mbox{E}_\kappa$}
-\Delta\varphi=p\lambda \alpha(x)(u^\kappa(x))^{p-1}\varphi\quad\mbox{in}\quad{\mathbb R}^N,
\quad
\varphi\in{\mathcal D}^{1,2}. 
\end{equation}
Then eigenvalue problem~\eqref{eq:E} has the first eigenvalue $\lambda^\kappa$ and the corresponding eigenfunction $\varphi^\kappa$ 
and the following properties hold:
\begin{eqnarray}
\notag
 & & \lambda^\kappa>0,\quad \varphi^\kappa(x)>0\quad\mbox{for a.a.~$x\in{\mathbb R}^N$},\\
\label{eq:5.2}
 & & \int_{{\mathbb R}^N}|\nabla\psi|^2\,dx\ge p\lambda^\kappa\int_{{\mathbb R}^N}\alpha (u^\kappa)^{p-1}\psi^2\,dx
 \quad\mbox{for $\psi\in{\mathcal D}^{1,2}$}.
\end{eqnarray}
\end{lemma}
Furthermore, we have the following lemmas.
\begin{Lemma}
\label{Lemma:5.2}
Assume the same conditions as in Theorem~{\rm\ref{Theorem:1.1}}. 
Let $\kappa\in{\mathcal K}_*$ and $\varphi^\kappa$ be as in the above. 
Then $\varphi^\kappa\in L^\infty_{0,-N+2}$ and 
\begin{equation}
\label{eq:5.3}
\varphi^\kappa(x)=p\lambda^\kappa \left[\Gamma*(\alpha (u^\kappa)^{p-1}\varphi^\kappa)\right](x)
\end{equation}
holds for a.a.~$x\in{\mathbb R}^N$.
Furthermore, there exists $C>0$ such that 
\begin{equation}
\label{eq:5.4}
\varphi^\kappa(x)\ge C(1+|x|)^{-N+2}\quad\mbox{for a.a.~$x\in{\mathbb R}^N$}. 
\end{equation}
\end{Lemma}
{\bf Proof.}
By \eqref{eq:5.1} we apply elliptic regularity theorems (see e.g. \cite{GT}*{Theorem~8.17}) to obtain $\varphi^\kappa\in L^\infty(\mathbb{R}^N)$. 
It follows from Lemma~\ref{Lemma:4.2} and $\varphi^\kappa\in{\mathcal D}^{1,2}$ that 
$(\varphi^{\kappa})^\sharp\in {\mathcal D}^{1,2}$ and 
$(\varphi^{\kappa})^\sharp$ satisfies
$$
-\Delta(\varphi^{\kappa})^\sharp=p\lambda^\kappa |x|^{-4}\left[\alpha (u^\kappa)^{p-1}\varphi^\kappa\right]^\sharp
=p\lambda^\kappa \beta ((u^{\kappa})^\sharp)^{p-1}(\varphi^{\kappa})^\sharp
\quad\mbox{in}\quad{\mathbb R}^N
$$
in the weak form. 
On the other hand, similarly to \eqref{eq:5.1}, 
by Lemma~\ref{Lemma:4.1} we see that
$$
\beta ((u^{\kappa})^\sharp)^{p-1}\in L^{q''}({\mathbb R}^N)\cap L^{\frac{N}{2}}({\mathbb R}^N)\quad\mbox{for some $q''>N/2$}.
$$
Then elliptic regularity theorems again yield
$(\varphi^{\kappa})^\sharp\in L^\infty$. 
This together with $\varphi^\kappa\in L^\infty$ implies that $\varphi^\kappa\in L^\infty_{0,-N+2}$.

We prove (\ref{eq:5.3}). Set
$$
\overline{\varphi}:=\lambda^\kappa\left[\Gamma*\left(p\alpha (u^\kappa)^{p-1}\varphi^\kappa\right)\right].
$$
By the same argument as in Lemma \ref{Lemma:3.3} 
we see that $\overline{\varphi}\in{\mathcal D}^{1,2}$ and 
$\overline{\varphi}$ satisfies
$-\Delta\overline{\varphi}=\lambda^\kappa\alpha p(u^\kappa)^{p-1}\varphi^\kappa$ in ${\mathbb R}^N$ in weak form. 
Then the uniqueness of weak solutions derives $\varphi^\kappa=\overline{\varphi}$.
Moreover, since $u^\kappa,\varphi^\kappa>0$ and $\alpha\not\equiv 0$ in ${\mathbb R}^N$, 
there exists $R>0$ such that $\alpha(u^\kappa)^{p-1}\varphi^\kappa\not\equiv 0$ in $B(0,R)$. 
This implies that
\begin{align*}
\varphi^\kappa(x) & 
\ge C\int_{B(0,R)}|x-y|^{-N+2}\alpha (u^\kappa)^{p-1}\varphi^\kappa\,dy\\
 & \ge C(|x|+R)^{-N+2}\int_{B(0,R)}\alpha(u^\kappa)^{p-1}\varphi^\kappa\,dy\ge C(1+|x|)^{-N+2}
\end{align*}
for a.a.~$x\in{\mathbb R}^N$. Thus \eqref{eq:5.4} holds, and Lemma~\ref{Lemma:5.2} follows. 
$\Box$

\begin{Lemma}
\label{Lemma:5.3}
Assume the same conditions as in Theorem~{\rm\ref{Theorem:1.1}}.
Then $\lambda^\kappa>1$ for all $\kappa\in(0,\kappa_*)$.
\end{Lemma}
{\bf Proof.}
Let $0<\kappa<\kappa'<\kappa_*$. 
Since 
\begin{equation}
\label{eq:5.5}
\mbox{the function $[0,\infty)\ni s\mapsto (t+s)^p-s^p$ is strictly increasing for any fixed $t>0$}, 
\end{equation}
by \eqref{eq:3.6} we have
\begin{align*}
(w^{\kappa'}+U^{\kappa'}_{j_*})^p-(U^{\kappa'}_{j_*-1})^p
 & =(w^{\kappa'}+V_{j_*}^{\kappa'}+U^{\kappa'}_{j_*-1})^p-(U^{\kappa'}_{j_*-1})^p\\
 & >(w^{\kappa'}+V_{j_*}^{\kappa'}+U^\kappa_{j_*-1})^p-(U^\kappa_{j_*-1})^p
\ge (w^{\kappa'}+U^\kappa_{j_*})^p-(U^\kappa_{j_*-1})^p.
\end{align*}
This together with \eqref{eq:3.7} implies that 
\begin{align*}
 & \lambda^\kappa\int_{{\mathbb R}^N}
p\alpha (u^\kappa)^{p-1}\varphi^\kappa(w^{\kappa'}-w^\kappa)\,dx
=\int_{{\mathbb R}^N} \nabla\varphi^\kappa\cdot\nabla(w^{\kappa'}-w^\kappa)\,dx\\
& \qquad\quad
=\int_{{\mathbb R}^N} \alpha\varphi^\kappa
\left\{(w^{\kappa'}+U^{\kappa'}_{j_*})^p-(U^{\kappa'}_{j_*-1})^p-(w^\kappa+U^\kappa_{j_*})^p+(U^\kappa_{j_*-1})^p\right\}\,dx\\
& \qquad\quad
>\int_{{\mathbb R}^N} \alpha\varphi^\kappa\left\{(w^{\kappa'}+U^\kappa_{j_*})^p-(w^\kappa+U^\kappa_{j_*})^p\right\}\,dx\\
& \qquad\quad
\ge\int_{{\mathbb R}^N} p\alpha (u^\kappa)^{p-1}\varphi^\kappa (w^{\kappa'}-w^\kappa)\,dx>0.
\end{align*}
This implies that $\lambda^\kappa>1$. 
Thus Lemma~\ref{Lemma:5.3} follows. 
$\Box$
\begin{Lemma}
\label{Lemma:5.4}
Assume the same conditions as in Theorem~{\rm\ref{Theorem:1.1}}.
Let $\kappa\in(0,\kappa^*)$. 
Then
$$
\int_{{\mathbb R}^N}p\alpha (u^\kappa)^{p-1}\psi^2\,dx
\le\int_{{\mathbb R}^N}|\nabla\psi|^2\,dx
\quad\mbox{for $\psi\in{\mathcal D}^{1,2}$}.
$$
\end{Lemma}
{\bf Proof.}
Let  $0<\kappa<\kappa'<\kappa^*$. 
Let $\{\eta_m\}_{m=1}^\infty$ be a sequence in $L^\infty_{\rm c}({\mathbb R}^N)\setminus\{0\}$ such that
\begin{equation}
\begin{split}
\label{eq:5.6}
 & 0\le\eta_1(x)\le\cdots\le\eta_m(x)\le\eta_{m+1}(x)\le\cdots\le p\alpha(x)u^\kappa(x)^{p-1},\\
 & \lim_{m\to\infty}\eta_m(x)=p\alpha(x)u^\kappa(x)^{p-1},
\end{split}
\end{equation}
for a.a.~$x\in{\mathbb R}^N$. 
Then, applying the same argument as in \cite{NS01}*{Lemma~B.2} with elliptic regularity theorems (see e.g. \cite{GT}*{Theorem~8.17}),
we find $\lambda_m>0$ and a function $\varphi_m\in {\mathcal D}^{1,2}\cap L^\infty({\mathbb R}^N)$ such that 
\begin{equation}
\label{eq:5.7}
-\Delta \varphi_m=\lambda_m \eta_m\varphi_m
\quad \mbox{in}\quad {\mathbb R}^N,
\qquad
\varphi_m>0
\quad \mbox{in}\quad {\mathbb R}^N.
\end{equation}
Set $W^\kappa_j:= U^\kappa_j-U^\kappa_{j_*}$ for $j\ge j_*+1$.
It follows from \eqref{eq:1.8} and Lemma~\ref{Lemma:2.4} that
$$
W^\kappa_j=\Gamma*\left(\alpha\left((U^\kappa_{j-1})^p-(U^\kappa_{j_*-1})^p\right)\right),
\qquad
W^\kappa_j=\sum_{k=j_*+1}^j V^\kappa_j\in L^\infty_{0,2-N}.
$$
Let $j\ge j_*+1$. Since
\[
W^\kappa_{j+1}\le\Gamma*\left(p\alpha(U^\kappa_j)^{p-1}(U^\kappa_j-U^\kappa_{j_*-1})\right)\le\Gamma*\left(p\alpha(U^\kappa_j)^{p-1}(W^\kappa_j+V^\kappa_{j_*})\right),
\]
applying the same argument as in the proof of Lemma~\ref{Lemma:3.3}, we see that $W^\kappa_j\in\mathcal{D}^{1,2}$ and
\begin{equation}\label{eq:5.8}
-\Delta W^{\kappa}_{j+1}=\alpha\left((U^\kappa_j)^p-(U^\kappa_{j_*-1})^p\right)\quad\mathrm{in}\quad\mathbb{R}^N
\end{equation}
in the weak form.
On the other hand, by \eqref{eq:5.5} we have
\[
(U^{\kappa'}_j)^p-(U^{\kappa'}_{j_*-1})^p
=(W^{\kappa'}_j+V_{j_*}^{\kappa'}+U^{\kappa'}_{j_*-1})^p-(U^{\kappa'}_{j_*-1})^p
> (W^{\kappa'}_j+V_{j_*}^{\kappa'}+U^\kappa_{j_*-1})^p-(U^\kappa_{j_*-1})^p.
\]
This together with \eqref{eq:5.7} and \eqref{eq:5.8} implies that
\begin{align*}
 & \lambda_m\int_{{\mathbb R}^N}\eta_m\varphi_m(W^{\kappa'}_{j+1}-W^\kappa_{j+1})\,dx
=\int_{{\mathbb R}^N}\nabla\varphi_m\cdot\nabla(W^{\kappa'}_{j+1}-W^\kappa_{j+1})\,dx\\
 & \qquad
 =\int_{{\mathbb R}^N}\alpha\varphi_m\left[(U^{\kappa'}_j)^p-(U^{\kappa'}_{j_*-1})^p-(U^\kappa_j)^p+(U^\kappa_{j_*-1})^p\right]\,dx 
\\
& \qquad
>\int_{{\mathbb R}^N}\alpha\varphi_m\left[(W^{\kappa'}_j+V^{\kappa'}_{j_*}+U^\kappa_{j_*})^p-(U^\kappa_j)^p\right]\,dx\\
& \qquad
>\int_{{\mathbb R}^N}\alpha\varphi_m\left[(W^{\kappa'}_j+U^\kappa_{j_*})^p-(W^\kappa_j+U^\kappa_{j_*})^p\right]\,dx
\ge\int_{{\mathbb R}^N}p\alpha(U^\kappa_j)^{p-1}\varphi_m (W^{\kappa'}_j-W^\kappa_j)\,dx
\end{align*}
for $j\ge j_*+1$. 
This together with \eqref{eq:5.6} implies that 
\begin{align*}
\infty & >\lambda_m\int_{{\mathbb R}^N}\eta_m\varphi_m(w^{\kappa'}-w^\kappa)\,dx
\ge \int_{{\mathbb R}^N}p\alpha(u^\kappa)^{p-1}\varphi_m(w^{\kappa'}-w^\kappa)\,dx\\
 & \ge \int_{{\mathbb R}^N}\eta_m\varphi_m(w^{\kappa'}-w^\kappa)\,dx>0.
\end{align*}
Then we observe that $\lambda_m\ge 1$, which together with \eqref{eq:5.7} yields
$$
\int_{{\mathbb R}^N}\eta_m\psi^2\,dx
\le\lambda_m\int_{{\mathbb R}^N}\eta_m\psi^2\,dx
\le\int_{{\mathbb R}^N}|\nabla\psi|^2\,dx,
\quad\psi\in{\mathcal D}^{1,2}.
$$
Therefore, by \eqref{eq:5.6} we obtain the desired conclusion. 
The proof is complete.
$\Box$
\vspace{5pt}

Now we are ready to complete the proof of Theorem~\ref{Theorem:1.1}.
\vspace{5pt}
\newline
{\bf Proof of Theorem~\ref{Theorem:1.1}.}
It follows from Lemma~\ref{Lemma:3.2} that $\kappa_*>0$. 
On the other hand, by Lemma~\ref{Lemma:5.4}, 
for any $\kappa\in(0,\kappa^*)$, we have
\[
0<\int_{{\mathbb R}^N}p\alpha (\kappa \Gamma*\mu)^{p-1}\psi^2\,dx
\le\int_{{\mathbb R}^N}p\alpha (u^\kappa)^{p-1}\psi^2\,dx\le\int_{{\mathbb R}^N}|\nabla\psi|^2\,dx,
\quad\psi\in{\mathcal D}^{1,2}\setminus\{0\}.
\] 
This implies that $\kappa^*<\infty$. Thus Theorem~\ref{Theorem:1.1} follows. 
$\Box$
\section{Uniform estimate for $w^\kappa$}
In this section we obtain a uniform $L^\infty(B(0,1))$ estimate of $\{w^\kappa\}_{\kappa\in(0,\kappa_*)}$ 
and $\{(w^\kappa)^\sharp\}_{\kappa\in(0,\kappa_*)}$. 
We first obtain a uniform energy estimates of $\{w^\kappa\}_{\kappa\in(0,\kappa_*)}$ in $B(0,3)$ 
to prove the following lemma. 
\begin{Lemma}
\label{Lemma:6.1}
Assume the same conditions as in Theorem~{\rm\ref{Theorem:1.1}}. 
Let $w^\kappa$ be as in \eqref{eq:3.4}. 
Then
\begin{equation}
\label{eq:6.1}
\sup_{0<\kappa<\kappa_*}\|w^\kappa\|_{L^{\frac{2N}{N-2}\nu}(B(0,2))}<\infty
\end{equation}
for all $\nu\ge 1$ with $\nu^2/(2\nu-1)<p$.
\end{Lemma}
{\bf Proof.} 
Let $\kappa\in(0,\kappa_*)$. 
Let $\nu\ge 1$ be such that $\nu^2/(2\nu-1)<p$.
Let $\eta\in C^\infty_c({\mathbb R}^N)$ be such that
$0\le\eta\le 1$ in ${\mathbb R}^N$, 
$\eta=1$ in $B(0,2)$, and 
$\eta=0$ outside $B(0,3)$. 
Let $\sigma>1$ be large enough such that
$$
1-\frac{1}{\sigma}\ge\frac{2\nu}{2\nu+p-1}. 
$$
Setting $\zeta=\eta^\sigma$, we have 
\begin{equation}
\label{eq:6.2}
\begin{split}
 & 0\le\zeta\le 1\quad\mbox{in}\quad{\mathbb R}^N,
\qquad
\zeta=1\quad\mbox{in}\quad B(0,2),
\qquad
\zeta=0\quad\mbox{outside}\quad B(0,3),\\
 & |\nabla\zeta|\le C\eta^{\sigma-1}=C\zeta^{1-\frac{1}{\sigma}}\le C\zeta^{\frac{2\nu}{2\nu+p-1}}\quad\mbox{in}\quad{\mathbb R}^N.
\end{split}
\end{equation}

Set $M:=\|V^{\kappa_*}_{j_*}\|_{L^\infty(B(0,3))}$ and $\overline{w}^\kappa:= w^\kappa+M$. 
Let $\varepsilon>0$ be small enough. 
Then we find $C_\varepsilon>0$ such that
\begin{equation}
\begin{aligned}
 & \int_{{\mathbb R}^N}|\nabla(\zeta (\overline{w}^\kappa)^\nu)|^2\,dx\\
 & =\int_{{\mathbb R}^N}\left(\nu^2\zeta^2(\overline{w}^\kappa)^{2\nu-2}
 |\nabla w^\kappa|^2+2\nu\zeta\nabla\zeta\cdot (\overline{w}^\kappa)^{2\nu-1}\nabla w^\kappa+|\nabla\zeta|^2(\overline{w}^\kappa)^{2\nu}\right)\,dx\\
&\le\int_{{\mathbb R}^N}\left(\frac{\nu^2(1+\varepsilon)}{2\nu-1}\nabla w^\kappa\cdot\nabla(\zeta^2(\overline{w}^\kappa)^{2\nu-1})
+C_\varepsilon|\nabla\zeta|^2(\overline{w}^\kappa)^{2\nu}\right)\,dx
\end{aligned}
\label{eq:6.3}
\end{equation}
for all $\kappa\in(0,\kappa_*)$. 
By \cite{IOS01}*{Lemma~5.1}, for any $\delta\in(0,1)$, we find $C_\delta>0$ such that 
$$
t^p-s^p\le (1+\varepsilon)t^{p-1}(t-s)+C_\delta s^{p-1+\delta}(t-s)^{1-\delta}
$$
for all $t$, $s\in[0,\infty)$ with $s\le t$.
Then, by \eqref{eq:3.7} we obtain
\begin{equation}
\begin{aligned}
 & \int_{{\mathbb R}^N}\nabla w^\kappa\cdot\nabla(\zeta^2(\overline{w}^\kappa)^{2\nu-1})\,dx\\
&=\int_{{\mathbb R}^N}\alpha\left((w^\kappa+U^{\kappa}_{j_*})^p-(U^\kappa_{j_*-1})^p\right)\zeta^2(\overline{w}^\kappa)^{2\nu-1}\,dx\\
&\le\int_{{\mathbb R}^N}\alpha\left((1+\varepsilon)(u^\kappa)^{p-1}(w^\kappa+V^\kappa_{j_*})
+C_\delta(U^\kappa_{j_*-1})^{p-1+\delta}(w^\kappa+V^\kappa_{j_*})^{1-\delta}\right)\zeta^2(\overline{w}^\kappa)^{2\nu-1}\,dx\\
&=\int_{{\mathbb R}^N}\alpha\left((1+\varepsilon)(u^\kappa)^{p-1}\zeta^2(\overline{w}^\kappa)^{2\nu}
+C_\delta(U^\kappa_{j_*-1})^{p-1+\delta}\zeta^2(\overline{w}^\kappa)^{2\nu-\delta}\right)\,dx.
\end{aligned}
\label{eq:6.4}\end{equation}
Furthermore, by Lemma \ref{Lemma:5.3} we see that
\begin{equation}
\label{eq:6.5}
\begin{split}
(1+\varepsilon)\int_{{\mathbb R}^N}\alpha (u^\kappa)^{p-1}\zeta^2(\overline{w}^\kappa)^{2\nu}\,dx
 & \le\frac{1+\varepsilon}{\lambda^\kappa p}\int_{{\mathbb R}^N}|\nabla(\zeta(\overline{w}^\kappa)^\nu)|^2\,dx\\
 & \le\frac{1+\varepsilon}{p}\int_{{\mathbb R}^N}|\nabla(\zeta(\overline{w}^\kappa)^\nu)|^2\,dx.
\end{split}
\end{equation}
By \eqref{eq:2.13} we find a small enough $\delta>0$ so that 
$$
q_\delta:=\frac{r}{p-1+\delta}>\frac{N}{2},\quad c_\delta:=(p-1+\delta)c_*+a>-2,\quad\mbox{and}\quad d_\delta:=(p-1+\delta)d_*+b<-2.
$$
Then, by Lemma~\ref{Lemma:2.1}, \eqref{eq:2.10}, and \eqref{eq:2.18} we obtain 
\begin{align*}
\sup_{0<\kappa<\kappa_*}\|\alpha(U^\kappa_{j_*})^{p-1+\delta}\|_{L^\frac{N}{2}(B(0,3))}
 & \le\|\alpha(U^{\kappa_*}_{j_*})^{p-1+\delta}\|_{L^\frac{N}{2}(B(0,3))}
\le C\|\alpha(U^{\kappa_*}_{j_*})^{p-1+\delta}\|_{L^{\frac{N}{2}}_{c_\delta,d_\delta}(B(0,3))}\\
 & \le C\|\alpha(U^{\kappa_*}_{j_*})^{p-1+\delta}\|_{L^{q_\delta}_{c_\delta,d_\delta}(B(0,3))}
 \le C\|U^{\kappa_*}_{j_*}\|_{L^r_{c_*,d_*}}^{p-1+\delta}
 \le C\kappa_*^{p-1+\delta}.
\end{align*}
This together with H\"{o}lder's inequality and the Sobolev inequality implies that
\begin{equation}
\label{eq:6.6}
\begin{split}
 & C_\delta\int_{{\mathbb R}^N}\alpha (U^\kappa_{j_*-1})^{p-1+\delta}\zeta^2(\overline{w}^\kappa)^{2\nu-\delta}\,dx\\
 & \le C_\delta C\|\zeta^2(\overline{w}^\kappa)^{2\nu-\delta}\|_{L^\frac{N}{N-2}({\mathbb R}^N)}
 \le C_\delta C\|\zeta(\overline{w}^\kappa)^\nu\|_{L^\frac{2N}{N-2}({\mathbb R}^N)}^{\frac{2\nu-\delta}{\nu}}\\
 & \le C_\delta C\|\nabla(\zeta(\overline{w}^\kappa)^\nu)\|_{L^2({\mathbb R}^N)}^{\frac{2\nu-\delta}{\nu}}
 \le \varepsilon \|\nabla(\zeta(\overline{w}^\kappa)^\nu)\|_{L^2({\mathbb R}^N)}^2+C
\end{split}
\end{equation}
for all $\kappa\in(0,\kappa_*)$. 
Combining \eqref{eq:6.4}, \eqref{eq:6.5}, and \eqref{eq:6.6}, we obtain
\begin{equation}
\int_{{\mathbb R}^N}\nabla w^\kappa\cdot\nabla(\zeta^2(\overline{w}^\kappa)^{2\nu-1})\,dx
\le\left(\frac{1+\varepsilon}{p}+\varepsilon\right)\int_{{\mathbb R}^N}|\nabla(\zeta(\overline{w}^\kappa)^\nu)|^2\,dx+C
\label{eq:6.7}\end{equation}
for all $\kappa\in(0,\kappa_*)$. 

On the other hand, 
since $\alpha$ is positive continuous in ${\mathbb R}^N\setminus\{0\}$, 
by \eqref{eq:6.2} we have 
\begin{equation}
\label{eq:6.8}
\begin{aligned}
 & \int_{{\mathbb R}^N}|\nabla\zeta|^2(\overline{w}^\kappa)^{2\nu}\,dx\\
 & =\int_{B(0,3)\setminus B(0,1)}|\nabla\zeta|^2(\overline{w}^\kappa)^{2\nu}\,dx
 \le C\int_{B(0,3)\setminus B(0,1)}\alpha^\frac{2\nu}{2\nu+p-1}|\nabla\zeta|^2(\overline{w}^\kappa)^{2\nu}\,dx\\
&\le C\int_{B(0,3)}\alpha^\frac{2\nu}{2\nu+p-1}\zeta^\frac{4\nu}{2\nu+p-1}(\overline{w}^\kappa)^{2\nu}\,dx
\le C\left(\int_{B(0,3)}\alpha\zeta^2(\overline{w}^\kappa)^{2\nu+p-1}\,dx\right)^\frac{2\nu}{2\nu+p-1}.
\end{aligned}
\end{equation}
Furthermore, we observe from \eqref{eq:6.5} and $\alpha\in L^{N/2}(B(0,3))$ that
\begin{equation}
\label{eq:6.9}
\begin{aligned}
& \int_{B(0,3)}\alpha((u^\kappa)^{p-1}+M^{p-1})\zeta^2(\overline{w}^\kappa)^{2\nu}\,dx\\
&\le C\int_{{\mathbb R}^N}|\nabla(\zeta(\overline{w}^\kappa)^\nu)|^2\,dx+\|\zeta^2(\overline{w}^\kappa)^{2\nu}\|_{L^\frac{N}{N-2}({\mathbb R}^N)}
\le C\int_{{\mathbb R}^N}|\nabla(\zeta(\overline{w}^\kappa)^\nu)|^2\,dx.
\end{aligned}
\end{equation}
By \eqref{eq:6.8} and \eqref{eq:6.9}, for any $\varepsilon'>0$, we see that
$$
\int_{{\mathbb R}^N}|\nabla\zeta|^2(\overline{w}^\kappa)^{2\nu}\,dx
\le\varepsilon'\int_{{\mathbb R}^N}|\nabla(\zeta(\overline{w}^\kappa)^\nu)|^2\,dx+C.
$$
Therefore, combining \eqref{eq:6.3}, \eqref{eq:6.7}, and \eqref{eq:6.8}, we obtain
$$
\int_{{\mathbb R}^N}|\nabla(\zeta (\overline{w}^\kappa)^\nu)|^2\,dx\\
\le\left[\frac{\nu^2(1+\varepsilon)}{2\nu-1}\left(\frac{1+\varepsilon}{p}+\varepsilon\right)+C_\varepsilon\varepsilon'\right]
 \int_{{\mathbb R}^N}|\nabla(\zeta(\overline{w}^\kappa)^\nu)|^2\,dx+C
$$
for all $\kappa\in(0,\kappa_*)$. 
Since $\nu^2/(2\nu-1)<p$, 
taking suitable small enough $\varepsilon$, $\varepsilon'>0$, we obtain 
$$
\sup_{0<\kappa<\kappa_*}\int_{{\mathbb R}^N}|\nabla (\zeta(w^\kappa)^\nu)|^2\,dx<\infty.
$$
This together with the Sobolev inequality implies \eqref{eq:6.1}. 
Thus Lemma~\ref{Lemma:6.1} follows.
$\Box$
\begin{Lemma}
\label{Lemma:6.2}
Assume that 
\begin{equation}
\label{eq:6.10}
\frac{\nu^2}{2\nu-1}<p<1+\frac{2(2+a_-)}{N-2}\nu\quad\mbox{for some $\nu>1$}.
\end{equation}
Then
\begin{equation}
\label{eq:6.11}
\sup_{0<\kappa<\kappa_*}\|w^\kappa\|_{L^\infty(B(0,1))}<\infty.
\end{equation}
\end{Lemma}
{\bf Proof.}
Let $\nu>1$ be as in \eqref{eq:6.10}. By \eqref{eq:1.1}, \eqref{eq:2.13}, and Lemma~\ref{Lemma:2.1}~(3) we have
$$
\alpha (U^{\kappa_*}_{j_*})^{p-1}\in L^q_{\overline{c},\overline{d}}\subset L^{q'}(\mathbb{R}^N)\quad\textrm{for some $q'>N/2$}.
$$
For any $\delta>0$, set
$$
\frac{1}{q_{\nu,\delta}}:=\max\left\{\frac{(N-2)(p-1)}{2N\nu}+\frac{-a_-+\delta}{N},\frac{1}{q'}\right\}.
$$
Taking small enough $\delta$, by \eqref{eq:6.10} we see that $q_{\nu,\delta}>N/2$. 
Then we observe from \eqref{eq:1.1} and \eqref{eq:6.1} that
$$
\sup_{0<\kappa<\kappa_*}\left\|\alpha (w^\kappa)^{p-1}\right\|_{L^{q_{\nu,\delta}}(B(0,2))}<\infty.
$$
Since $\alpha (U^{\kappa_*}_{j_*})^{p-1}\in L^{q'}(B(0,2))$ 
and $q_{\nu,\delta}\le q'$, we deduce that
\begin{align*}
 & \sup_{0<\kappa<\kappa_*}\|\alpha (u^{\kappa})^{p-1}\|_{L^{q_{\nu,\delta}}(B(0,2))}\\
 & \le C\sup_{0<\kappa<\kappa_*}\left\|\alpha (w^\kappa)^{p-1}\right\|_{L^{q_{\nu,\delta}(B(0,2))}}+C\|\alpha (U^{\kappa_*}_{j_*})^{p-1}\|_{L^{q'}(B(0,2))}
<\infty.
\end{align*}
Since
$$
0\le \alpha((w^\kappa+U^\kappa_{j_*})^p-(U^\kappa_{j_*-1})^p)
\le p\alpha (u^\kappa)^{p-1}w^\kappa
\quad\mbox{in}\quad{\mathbb R}^N,
$$
we apply elliptic regularity theorems (see e.g. \cite{GT}*{Theorem~8.17}) to elliptic equation~\eqref{eq:3.7} to obtain \eqref{eq:6.11}. 
Thus Lemma~\ref{Lemma:6.2} follows. 
$\Box$\vspace{5pt}

At the end of this section, 
combing Lemma~\ref{Lemma:6.2} with the Kelvin transform, 
we obtain a decay estimate of $w^\kappa$ at the space infinity. 
\begin{lemma}
\label{Lemma:6.3}
Let $a^\sharp$ be as in Lemma~{\rm\ref{Lemma:4.1}}, that is, $a^\sharp:=(N-2)(p-1)-4-b$.
Assume that 
\begin{equation}
\label{eq:6.12}
\frac{\nu^2}{2\nu-1}<p<1+\frac{2(2+(a^\sharp)_-)}{N-2}\nu\quad\mbox{for some $\nu>1$}.
\end{equation}
Then
\begin{equation}
\label{eq:6.13}
\sup_{0<\kappa<\kappa_*}\sup_{|x|>1}|x|^{N-2}|w^\kappa(x)|<\infty.
\end{equation}
\end{lemma}
{\bf Proof.}
By Lemma~\ref{Lemma:4.1}, \eqref{eq:4.1}, and \eqref{eq:6.12} 
we apply Lemma~\ref{Lemma:6.2} to integral equation~\eqref{eq:I} to obtain 
$$
\sup_{0<\kappa<\kappa_*}\|(w^\kappa)^\sharp\|_{L^\infty(B(0,1))}<\infty.
$$
This implies \eqref{eq:6.13}, and Lemma~\ref{Lemma:6.3} follows. 
$\Box$
\section{Proof of Theorem~\ref{Theorem:1.2}}
For the proof of Theorem~\ref{Theorem:1.2}, 
we prepare the following lemma on relations~\eqref{eq:6.10} and \eqref{eq:6.12}.
\begin{Lemma}
\label{Lemma:7.1} 
Assume that $a>-2$ and $b\in{\mathbb R}$.
\begin{itemize}
  \item[{\rm (i)}] 
  Relation~\eqref{eq:6.10} holds if and only if $1<p<p^*(a_-)$.
  \item[{\rm (ii)}]
  Relation~\eqref{eq:6.12} holds if and only if $p_*(b)<p<p^*(0)$.  
\end{itemize}
\end{Lemma}
{\bf Proof of assertion~(i).}
We first show that 
\begin{equation}
\label{eq:7.1}
I:=\left\{\nu>1\,:\,\frac{\nu^2}{2\nu-1}<1+\frac{2(2+a_-)}{N-2}\nu\right\}=(1,\nu^*),
\end{equation}
where
$$
\nu^*=\left\{\begin{aligned}
&\frac{N-2}{N-a_--4-\sqrt{(a_-+2)(2N+a_--2)}}& &\mathrm{if}\  N>10+4a_-,\\
&\infty& &\mathrm{otherwise}.
\end{aligned}\right.
$$
We notice that, 
for any $\nu>1$, $\nu\in I$ is equivalent to
\begin{equation}
\label{eq:7.2}
h(\nu):=(N-4a_--10)\nu^2-2(N-a_--4)\nu+(N-2)<0.
\end{equation}
It follows from $a>-2$ that 
\begin{equation}
\label{eq:7.3}
\begin{split}
 & D:=(N-a_--4)^2-(N-2)(N-4a_--10)=(a_-+2)(2N+a_--2)>0,\\
 & h(1)=-2a_--4<0. 
\end{split}
\end{equation}

We consider the case of $N>10+4a_-$. 
By \eqref{eq:7.2} and \eqref{eq:7.3} we see that 
$\nu\in I$ is equivalent to 
\begin{equation*}
\begin{split}
\nu<\frac{N-a_--4+\sqrt{D}}{N-4a_--10}
 & =\frac{(N-a_--4)^2-D}{(N-4a_--10)(N-a_--4-\sqrt{D})}\\
 & =\frac{(N-2)(N-4a_--10)}{(N-4a_--10)(N-a_--4-\sqrt{D})}
=\nu^*.
\end{split}
\end{equation*}
Thus \eqref{eq:7.1} holds and $h(\nu^*)=0$ in the case of $N>10+4a_-$. 

Next, we consider the case of $N<10+4a_-$. 
Since
$$
\frac{N-a_--4}{N-4a_--10}=1+\frac{3a_-+6}{N-4a_--10}<1,
$$
by \eqref{eq:7.2} and \eqref{eq:7.3} we see that $\nu\in I$ for all $\nu>1$. 
Thus \eqref{eq:7.1} holds in the case of $N<10+4a_-$.
In the case of $N=10+4a_-$, 
since
$$
N-a_--4=N-4a_-+10+3a_-+6=3a_-+6>0, 
$$
by \eqref{eq:7.2} and \eqref{eq:7.3} we see that $\nu\in I$ for all $\nu>1$. 
Thus \eqref{eq:7.1} holds in the case of $N=10+4a_-$, 
and we obtain \eqref{eq:7.1}.

We complete the proof of assertion~(i). 
If $N\le 10+4a_-$, then $\nu^*=\infty$ and $p^*(a_-)=\infty$. 
This together with \eqref{eq:7.1} implies assertion~(i) in the case of $N\le 10+4a_-$. 
If $N>10+4a_-$, then $h(\nu^*)=0$, which yields
$$
\frac{(\nu^*)^2}{2\nu^*-1}=1+\frac{2(2+a_-)}{N-2}\nu^*
=1+\frac{2(2+a_-)}{N-a_--4-\sqrt{(a_-+2)(2N+a_--2)}}
=p^*(a_-).
$$
Since 
\begin{equation}
\label{eq:7.4}
\mbox{the function $(1,\infty)\ni\nu\mapsto \displaystyle{\frac{\nu^2}{2\nu-1}}$ is monotone increasing},
\end{equation} 
by \eqref{eq:7.1} we see that 
\begin{equation}
\label{eq:7.5}
\begin{split}
 & \inf_{\nu\in I}\frac{\nu^2}{2\nu-1}=1,\\
 & \sup_{\nu\in I}\frac{\nu^2}{2\nu-1}=\frac{(\nu^*)^2}{2\nu^*-1}=p^*(a_-),\\
 & 1+\sup_{\nu\in I}\frac{2(2+a_-)}{N-2}\nu=1+\frac{2(2+a_-)}{N-2}\nu^*=p^*(a_-).
\end{split}
\end{equation}
Therefore, if $p$ satisfies relation~\eqref{eq:6.10}, 
then there exists $\nu\in I$ such that
$$
1\le \frac{\nu^2}{2\nu-1}<p<1+\frac{2(2+a_-)}{N-2}\nu\le p^*(a_-).
$$

Conversely, assume $1<p<p^*(a_-)$. 
By \eqref{eq:7.5} we apply the intermediate value theorem to obtain $\nu_p\in I$ such that 
$$
\frac{\nu_p^2}{2\nu_p-1}=p.
$$
Since $\nu_p\in I$, 
we find a small enough $\delta>0$ such that
$$
\nu_p-\delta>1,\quad \frac{\nu_p^2}{2\nu_p-1}<1+\frac{2(2+a_-)}{N-2}(\nu_p-\delta).
$$
Therefore, by \eqref{eq:7.4} we have
\[
\frac{(\nu_p-\delta)^2}{2(\nu_p-\delta)-1}<\frac{\nu_p^2}{2\nu_p-1}=p< 1+\frac{2(2+a_-)}{N-2}(\nu_p-\delta).
\]
Thus relation~\eqref{eq:6.10} holds with $\nu=\nu_p-\delta$. 
Then assertion~(i) holds in the case of $N>10+4a_-$, and the proof of assertion~(i) is complete.
$\Box$\vspace{3pt}
\newline
{\bf Proof of assertion~(ii) in the case of $a^\sharp\ge 0$.}
Let $a^\sharp\ge 0$, that is, 
$$
p\ge\frac{N+2+b}{N-2}.
$$
By assertion~(i) we see that 
relation~\eqref{eq:6.12} holds if and only if $1<p<p^*(0)$. 
For the proof of assertion~(ii), 
it suffices to prove 
\begin{equation}
\label{eq:7.6}
\left[\frac{N+2+b}{N-2},\infty\right)\cap(1,p^*(0))
=\left[\frac{N+2+b}{N-2},\infty\right)\cap(p_*(b),p^*(0)).
\end{equation}
Notice that $p_*(b)\ge 1$ (see \eqref{eq:1.5}). 
If $b\le -2$, then $p_*(b)=1$, and \eqref{eq:7.6} holds. 

Consider the case of $b>-2$. Since $p_*(b)\ge 1$, if 
$$
p_*(b)<\frac{N+2+b}{N-2},
$$
then \eqref{eq:7.6} holds.
Furthermore, 
\begin{align*}
\frac{N+2+b}{N-2}-p_*(b)
& =\frac{N+2+b}{N-2}-\left(1+\frac{2(2+b)}{N-b-4+\sqrt{D'}}\right)\\
&=\frac{(b+4)(N-b-4+\sqrt{D'})-2(N-2)(2+b)}{(N-2)(N-b-4+\sqrt{D'})}\\
&=\frac{(b+4)\sqrt{D'}-b^2-(N+4)b-8}{(N-2)(N-b-4+\sqrt{D'})},
\end{align*}
where $D':=(N+b)^2-(N-2)^2=(b+2)(2N+b-2)$.
These mean that \eqref{eq:7.6} holds if either 
$$
({\rm B1})\quad E:=(b+4)^2D'-(b^2+(N+4)b+8)^2> 0
\qquad\mbox{or}\qquad
({\rm B2})\quad b^2+(N+4)b+8< 0.
$$
Here
\begin{equation}
\label{eq:7.7}
\begin{aligned}
 E &=(b+4)^2(b+2)(b+2N-2)-(b^2+(N+4)b+8)^2\\
&=(-N^2+12N-20)b^2+(48N-96)b+64N-128\\
&=(N-2)((10-N)b^2+48b+64).
\end{aligned}
\end{equation}
Let 
\begin{equation}
\label{eq:7.8}
b_\pm:=\frac{-24\pm8\sqrt{N-1}}{10-N}=-\frac{8}{3\pm\sqrt{N-1}},
\end{equation}
which are the roots of the algebra equation $(10-N)x^2+48x+64=0$. 
Furthermore, let 
$$
c_\pm:=\frac{-N-4\pm\sqrt{N^2+8N-16}}{2}=-\frac{16}{N+4\pm\sqrt{N^2+8N-16}},
$$
which are also the roots of the algebra equation $x^2+(N+4)x+8=0$. 
Since $c_-<-(N+4)/2<-2$ and $b>-2$, we see that
\begin{equation}
\label{eq:7.9}
 \mbox{(B2) holds if and only if  
 $b\in\left(-2,c_+\right)$}.
\end{equation}
Furthermore, since
$N+4>6$ and $N^2+8N-16=(N+6)(N-2)+4N-4>4N-4$,
we have
\begin{equation}
\label{eq:7.10}
b_+=-\frac{16}{6+2\sqrt{N-1}}<-\frac{16}{N+4+\sqrt{N^2+8N-16}}=c_+.
\end{equation}
We divide the proof of  \eqref{eq:7.6} into three cases $N\le 9$, $N=10$, and $N\ge 11$.
\vspace{2pt}
\newline
\underline{Case of $N\le 9$} 
By \eqref{eq:7.7} and \eqref{eq:7.8} we see that, 
if $b\not\in[b_-,b_+]$, then $E\ge 0$, that is, (B1) holds. 
On the other hand, 
if $b\in[b_-,b_+]$, then, by \eqref{eq:7.10} we have $b\in(-2,c_+)$, which together with \eqref{eq:7.9} implies that (B2) holds.
Thus \eqref{eq:7.6} holds in the case of $N\ge 9$.
\vspace{2pt}
\newline
\underline{Case of $N=10$} 
By \eqref{eq:7.7} we see that, if $b> -4/3$, then $E> 0$, that is, (B1) holds.  
On the other hand, since 
$$
-\frac{4}{3}=-\frac{8}{3+\sqrt{9}}=b_+<c_+,
$$
if $b\le-4/3$, then $b\in(-2,c_+)$, which together with \eqref{eq:7.9} implies that (B2) holds.
Thus \eqref{eq:7.6} holds in the case of $N=10$.
\vspace{2pt}
\newline
\underline{Case of $N\ge 11$} 
By \eqref{eq:7.7} we see that, if $b\in (b_+.b_-)$, then $E> 0$, that is, (B1) holds. 
On the other hand, 
if $b\le b_+$, then, by \eqref{eq:7.10} we have $b\in(-2,c_+)$, which together with \eqref{eq:7.9} implies that (B2) holds. 
Thus \eqref{eq:7.6} holds in the case of $N\ge 11$ with $b<b_-$.

Consider the case of $b\ge b_-$. 
Then, for any $p>1$ with
$$
p\ge\frac{N+2+b}{N-2}=1+\frac{4+b}{N-2},
$$
we have
\begin{align*}
p& \ge1+\frac{4+b_-}{N-2}=1+\frac{4(N-10)+24+8\sqrt{N-1}}{(N-2)(N-10)}\\
&=1+\frac{4N-16+8\sqrt{N-1}}{(N-2)(N-10)}=1+\frac{4(N-4+2\sqrt{N-1})(N-4-2\sqrt{N-1})}{(N-2)(N-10)(N-4-2\sqrt{N-1})}\\
&=1+\frac{4((N-4)^2-4(N-1))}{(N-2)(N-10)(N-4-2\sqrt{N-1})}\\
&=1+\frac{4}{N-4-2\sqrt{N-1}}=p^*(0).
\end{align*}
This implies that 
$$
\left[\frac{N+2+b}{N-2},\infty\right)\cap(1,p^*(0))
=\left[\frac{N+2+b}{N-2},\infty\right)\cap(p_*(b),p^*(0))=\emptyset.
$$
Thus \eqref{eq:7.6} holds in the case of $N\ge 11$ with $b>b_-$.
Therefore assertion~(ii) follows in the case of $a^\sharp\ge 0$.
$\Box$\vspace{3pt}
\newline
{\bf Proof of assertion~(ii) in the case of $a^\sharp<0$.}
Let $a^\sharp<0$.  
Then relation~\eqref{eq:6.12} holds if and only if 
$$
p>\max\left\{\frac{\nu^2}{2\nu-1},1+\frac{2(2+b)\nu}{(N-2)(2\nu-1)}\right\}\quad\mbox{for some $\nu>1$}.
$$
Thus it suffices to prove that
\begin{equation}
\label{eq:7.11}
\inf_{\nu>1}\max\left\{\frac{\nu^2}{2\nu-1},1+\frac{2(2+b)\nu}{(N-2)(2\nu-1)}\right\}=p_*(b).
\end{equation}
Consider the case of $b\le -2$. By \eqref{eq:7.4} we have
$$
\frac{\nu^2}{2\nu-1}> 1\ge1+\frac{2(2+b)\nu}{(N-2)(2\nu-1)}
\quad\mbox{for all $\nu>1$},
$$
which implies that 
$$
\inf_{\nu>1}\max\left\{\frac{\nu^2}{2\nu-1},1+\frac{2(2+b)\nu}{(N-2)(2\nu-1)}\right\}=\inf_{\nu>1}\frac{\nu^2}{2\nu-1}=1=p_*(b).
$$
Thus \eqref{eq:7.11} holds in the case of $b\le -2$. 

Consider the case of $b>-2$. 
Since the function 
$$
(1,\infty)\ni\nu\mapsto 1+\frac{2(2+b)}{(N-2)(2\nu-1)}
$$
is strictly monotone decreasing, by \eqref{eq:7.4} we apply the intermediate value theorem 
to find a unique  $\nu_*>1$ such that
$$
\frac{\nu_*^2}{2\nu_*-1}=1+\frac{2(2+b)\nu_*}{(N-2)(2\nu_*-1)},
$$
that is, 
$$
\nu_*=\frac{N+b+\sqrt{D'}}{N-2}.
$$
Here we used that 
$$
\frac{N+b-\sqrt{D'}}{N-2}<1.
$$
Since 
\begin{align*}
\frac{\nu_*}{(N-2)(2\nu_*-1)} & =\frac{N+b+\sqrt{D'}}{(N-2)(N+2b+2+2\sqrt{D'})}\\
 & =\frac{(N+b+\sqrt{D'})(N+b-\sqrt{D'})}{(N-2)(N+b+\sqrt{D'}+b+2+\sqrt{D'})(N+b-\sqrt{D'})}\\
 & =\frac{N-2}{(N+b)^2-D'+(b+2+\sqrt{D'})(N+b-\sqrt{D'})}\\
 & =\frac{N-2}{(N-2)^2+(b+2)(N+b)+(N-2)\sqrt{D'}-D'}\\
 & =\frac{N-2}{2(N-2)^2-(N+b)(N-2)+(N-2)\sqrt{D'}}\\
 & =\frac{1}{2(N-2)-(N+b)+\sqrt{D'}}\\
 & =\frac{1}{N-b-4+\sqrt{(b+2)(2N+b-2)}},
\end{align*}
we deduce that
$$
\inf_{\nu>1}\max\left\{\frac{\nu^2}{2\nu-1},1+\frac{2(2+b)\nu}{(N-2)(2\nu-1)}\right\}=1+\frac{2(2+b)\nu_*}{(N-2)(2\nu_*-1)}=p_*(b).
$$
Thus \eqref{eq:7.11} holds in the case of $b>-2$, 
and we obtain assertion~(ii) in the case of $a^\sharp<0$. 
Therefore the proof of Lemma~\ref{Lemma:7.1} is complete.
$\Box$
\begin{Lemma}
\label{Lemma:7.2}
Assume the same conditions as in Theorem~{\rm\ref{Theorem:1.2}}. 
Then $\kappa_*\in{\mathcal K}_*$ and $\lambda^{\kappa_*}\ge 1$. 
\end{Lemma}
{\bf Proof.}
By Lemmas~\ref{Lemma:6.2} and \ref{Lemma:6.3} we see that
$w^\kappa(x)\le C\omega_{0,-N+2}(x)$ for a.a.~$x\in{\mathbb R}^N$ and all~$\kappa\in(0,\kappa_*)$. 
This together with \eqref{eq:3.6} implies that the limit 
$$
w(x):=\lim_{\kappa\to\kappa_*-0}w^\kappa(x)
$$
exists for a.a.~$x\in{\mathbb R}^N$ and $w\in L^\infty_{0,-N+2}$. 
On the other hand, 
it follows from \eqref{eq:2.10} that $U^\kappa_j\to U^{\kappa_*}_j$ as $\kappa\to\kappa_*$ for $j=0,1,2,\dots$. 
Then we see that $w$ satisfies \eqref{eq:3.5} with $\kappa=\kappa_*$, 
that is, the function $u=w+U^{\kappa_*}_{j_*}$ is a solution of \eqref{eq:P} with $\kappa=\kappa_*$ 
and $u\in L^r_{c_*,d_*}$. 
Thus $\kappa_*\in{\mathcal K}_*$. 
Furthermore, by Lemma~\ref{Lemma:5.2} and \eqref{eq:5.2} we have
$$
\int_{{\mathbb R}^N}|\nabla\psi|^2\,dx>p\int_{{\mathbb R}^N}\alpha p(u^\kappa)^{p-1}\psi^2\,dx
\ge p\int_{{\mathbb R}^N}\alpha p(U^\kappa_j)^{p-1}\psi^2\,dx
$$
for all $\psi\in{\mathcal D}^{1,2}$, $\kappa\in(0,\kappa_*)$, and $j=0,1,2,\dots$. 
Then we obtain
$$
\int_{{\mathbb R}^N}|\nabla\psi|^2\,dx
\ge p\int_{{\mathbb R}^N}\alpha p(U^{\kappa_*}_j)^{p-1}\psi^2\,dx
$$
for all $\psi\in{\mathcal D}^{1,2}$ and $j=0,1,2,\dots$. 
This together with Lemma~\ref{Lemma:3.1} implies that 
$$
\int_{{\mathbb R}^N}|\nabla\psi|^2\,dx
\ge p\int_{{\mathbb R}^N}\alpha p(u^{\kappa_*})^{p-1}\psi^2\,dx
$$
for all $\psi\in{\mathcal D}^{1,2}$.
Thus $\lambda^{\kappa_*}\ge 1$, and Lemma~\ref{Lemma:7.2} follows.
$\Box$\vspace{5pt}

Next, we characterize $\kappa_*$ using the eigenvalues $\{\lambda^\kappa\}$.
\begin{Lemma}
\label{Lemma:7.3}
Assume the same conditions as in Theorem~{\rm\ref{Theorem:1.2}}. 
If $0<\kappa\le\kappa_*$ and $\lambda^\kappa>1$, then $\kappa<\kappa_*$. 
Furthermore, if $\kappa=\kappa_*$, then $\lambda^{\kappa_*}=1$.
\end{Lemma}
{\bf Proof.}
Let $\delta\in(0,1)$ and $k>0$ be such that  $\varepsilon:=k\delta\in(0,1)$. 
Let $0<\kappa\le\kappa_*$ and assume that $\lambda^\kappa>1$. 
We show that the function
$$
\overline{u}:= u^\kappa-U^{\kappa}_{j_*}+U^{(1+\varepsilon)\kappa}_{j_*}+\delta\varphi^\kappa
$$ 
is a supersolution to problem~\eqref{eq:P} with $\kappa$ replaced by $(1+\varepsilon)\kappa$.
It follows that
\begin{equation}
\label{eq:7.12}
\begin{aligned}
 & \overline{u}-\Gamma*(\alpha\overline{u}^p)-(1+\varepsilon)\kappa\Gamma*\mu\\ 
 & =u^\kappa-U^{\kappa}_{j_*}-\Gamma*\left[\alpha\left(\overline{u}^p-\left(U^{(1+\varepsilon)\kappa}_{j_*-1}\right)^p\right)\right]+\delta\varphi^\kappa\\
 &=-\Gamma*
 \left[\alpha\left(\overline{u}^p-(u^\kappa)^p-\left(U^{(1+\varepsilon)\kappa}_{j_*-1}\right)^p+\left(U^\kappa_{j_*-1}\right)^p\right)\right]+\delta\varphi^\kappa.
\end{aligned}
\end{equation}
By \eqref{eq:2.12}, \eqref{eq:2.18}, \eqref{eq:3.1}, and \eqref{eq:5.4} 
we have
\begin{equation}
\label{eq:7.13}
\begin{split}
\overline{u}-\left(u^\kappa-U^\kappa_{j_*-1}+U^{(1+\varepsilon)\kappa}_{j_*-1}\right)
 & =V^{(1+\varepsilon)\kappa}_{j_*}-V^\kappa_{j_*}+\delta\varphi^\kappa>0,\\
\overline{u}^p-\left(u^\kappa-U^\kappa_{j_*-1}+U^{(1+\varepsilon)\kappa}_{j_*-1}\right)^p&\le p\overline{u}^{p-1}(V^{(1+\varepsilon)\kappa}_{j_*}-V^\kappa_{j_*}+\delta\varphi^\kappa)\\
&\le p(u^\kappa+C\varepsilon U^\kappa_{j_*}+\delta\varphi^\kappa)^{p-1}(C\varepsilon V^\kappa_{j_*}+\delta\varphi^\kappa)\\
&\le (1+C(\delta+\varepsilon))p(u^\kappa)^{p-1}(C\varepsilon V^\kappa_{j_*}+\delta\varphi^\kappa)\\
&\le (1+C(\delta+\varepsilon))p(u^\kappa)^{p-1}(C\varepsilon\omega_{0,-N+2}+\delta\varphi^\kappa)\\
&\le (1+C(\delta+\varepsilon))(\delta+C\varepsilon)p(u^\kappa)^{p-1}\varphi^\kappa\\
&=(1+C(1+k)\delta)(1+Ck)\delta p(u^\kappa)^{p-1}\varphi^\kappa.
\end{split}
\end{equation}
On the other hand, it follows from $p>1$ that 
\begin{equation}
\label{eq:7.14}
\begin{split}
 & (s+t+v)^p-(t+v)^p-(s+v)^p+v^p=\int_0^s p((\sigma+t+v)^{p-1}-(\sigma+v)^{p-1})\,d\sigma\\
 & \qquad\quad
 =\int_0^t\int_0^s p(p-1)(\sigma+\tau+v)^{p-2}\,d\sigma\,d\tau\\
 & \qquad\quad
\le\begin{cases}
p(p-1)(s+t+v)^{p-2}st&\quad\mathrm{if}\quad p\ge 2,\vspace{3pt}\\
p(p-1)v^{p-2}st&\quad\mathrm{if}\quad 1<p<2
\end{cases}\\
 & \qquad\quad
 =p(p-1)\max\{(s+t+v)^{p-2},v^{p-2}\}st
\end{split}
\end{equation}
for $s,t,v>0$.
Applying \eqref{eq:7.14} with 
$$
t=u^\kappa-U^\kappa_{j_*-1},\quad
s=U^{(1+\varepsilon)\kappa}_{j_*-1}-U^\kappa_{j_*-1},\quad
v=U^\kappa_{j_*-1},
$$
by \eqref{eq:2.12}, \eqref{eq:2.18}, and \eqref{eq:5.4} we have
\begin{equation}
\label{eq:7.15}
\begin{aligned}
 & \left(u^\kappa-U^\kappa_{j_*-1}+U^{(1+\varepsilon)\kappa}_{j_*-1}\right)^p-(u^\kappa)^p
-\left(U^{(1+\varepsilon)\kappa}_{j_*-1}\right)^p+\left(U^\kappa_{j_*-1}\right)^p\\
 & \le p(p-1)\max\left\{\left(u^\kappa-U^{\kappa}_{j_*-1}+U^{(1+\varepsilon)\kappa}_{j_*-1}\right)^{p-2},\left(U^\kappa_{j_*-1}\right)^{p-2}\right\}\\
 & \qquad\qquad\qquad\qquad\qquad\qquad\qquad\qquad
 \times\left(U^{(1+\varepsilon)\kappa}_{j_*-1}-U^\kappa_{j_*-1}\right)(u^\kappa-U^\kappa_{j_*-1})\\
 & \le C\varepsilon\max\left\{\left(u^\kappa-U^\kappa_{j_*-1}+U^{(1+\varepsilon)\kappa}_{j_*-1}\right)^{p-2},\left(U^\kappa_{j_*-1}\right)^{p-2}\right\}
 U^\kappa_{j_*-1}(w^\kappa+V_{j_*}^\kappa)\\
 & \le C\varepsilon\max\left\{\left(u^\kappa-U^\kappa_{j_*-1}+U^{(1+\varepsilon)\kappa}_{j_*-1}\right)^{p-1},\left(U^\kappa_{j_*-1}\right)^{p-1}\right\}
(w^\kappa+V_{j_*}^\kappa)\\
 & = C\varepsilon\left(u^\kappa-U^\kappa_{j_*-1}+U^{(1+\varepsilon)\kappa}_{j_*-1}\right)^{p-1}(w^\kappa+V_{j_*}^\kappa)\\
 & \le C\varepsilon(u^\kappa+C\varepsilon U^\kappa_{j_*-1})^{p-1}\omega_{0,-N+2}\le C\varepsilon (u^\kappa)^{p-1}\varphi^\kappa.
\end{aligned}
\end{equation}
Combing \eqref{eq:7.13} and \eqref{eq:7.15}, we obtain
$$
\overline{u}^p-(u^\kappa)^p-\left(U^{(1+\varepsilon)\kappa}_{j_*-1}\right)^p+\left(U^\kappa_{j_*-1}\right)^p\le ((1+C(1+k)\delta)(1+Ck)+Ck)\delta(u^\kappa)^{p-1}\varphi^\kappa.
$$
This together with \eqref{eq:5.3} and \eqref{eq:7.12} implies that
\begin{align*}
 & \overline{u}-\Gamma*\alpha\overline{u}^p-(1+\varepsilon)\kappa\Gamma*\mu\\
 &\ge -((1+C(1+k)\delta)(1+Ck)+Ck)\delta\left[\Gamma*(p \alpha(u^\kappa)^{p-1}\varphi^\kappa)\right]+\delta\varphi^\kappa\\
&=\left(1-((1+C(1+k)\delta)(1+Ck)+Ck)(\lambda^\kappa)^{-1}\right)\delta\varphi^\kappa.
\end{align*}
Since $\lambda^\kappa>1$, 
taking small enough $\delta$ and $k$, we see that
$$
\overline{u}-\Gamma*(\alpha\overline{u}^p)-(1+\varepsilon)\kappa\Gamma*\mu> 0\quad\mbox{in}\quad{\mathbb R}^N.
$$
Therefore we see that $\overline{u}$ is a supersolution to problem~\eqref{eq:P} 
with $\kappa$ replaced by $(1+\varepsilon)\kappa$. 
Since $0<\kappa\le\kappa_*$, 
combining Lemma~\ref{Lemma:2.4}, Lemma~\ref{Lemma:7.2}, and \eqref{eq:1.12},
we have $\overline{u}\in L^r_{c_*,d_*}$. 
Then Lemma~\ref{Lemma:3.1} implies that $(1+\varepsilon)\kappa\in{\mathcal K}_*$, 
and we deduce that $\kappa<(1+\varepsilon)\kappa\le\kappa_*$. 
Furthermore, thanks to Lemma~\ref{Lemma:7.2}, 
we observe that $\lambda^{\kappa_*}=1$. 
Thus Lemma~\ref{Lemma:7.3} follows.
$\Box$
\vspace{5pt}

Now we are ready to complete the proof of Theorem~\ref{Theorem:1.2}.
\vspace{3pt}
\newline
{\bf Proof  of Theorem~\ref{Theorem:1.2}.}
We prove the uniqueness of solutions to problem~\eqref{eq:P} with $\kappa=\kappa_*$. 
Let $\tilde{u}$ be a solution to problem~\eqref{eq:P} with $\kappa=\kappa_*$. 
Since $u^{\kappa_*}$ is a minimal solution, we see that $\tilde{u}\ge u^{\kappa_*}$. 
Let $\{z_\ell\}_{\ell=0}^\infty$ be a sequence in $L^\infty_c({\mathbb R}^N\setminus\{0\})$ such that
$$
0\le z_0(x)\le \cdots\le z_\ell(x)\le z_{\ell+1}(x)\le\cdots\le\tilde{u}(x)-u^{\kappa_*}(x),
\quad
\lim_{\ell\to\infty}z_\ell(x)=\tilde{u}(x)-u^{\kappa_*}(x),
$$
for a.a.~$x\in{\mathbb R}^N$. 
Set
$$
Z_\ell:=\Gamma*\left[\alpha\left((u^{\kappa_*}+z_\ell)^p-(u^{\kappa_*})^p\right)\right].
$$
Since
$$
0\le \alpha((u^{\kappa_*}+z_\ell)^p-(u^{\kappa_*})^p)\le \alpha p(u^{\kappa_*}+z_\ell)^{p-1}z_\ell\in L^{q'}_c({\mathbb R}^N\backslash\{0\})
$$
for some $q'>N/2$
(see \eqref{eq:5.1}), 
by the same argument as Lemma~\ref{Lemma:3.3}
we obtain 
$Z_\ell\in L^\infty_{0,-N+2}$,
$Z_\ell\in{\mathcal D}^{1,2}$, and
\begin{equation}
-\Delta Z_\ell=\alpha((u^{\kappa_*}+z_\ell)^p-(u^{\kappa_*})^p)\quad \mathrm{in}\ {\mathbb R}^N
\label{eq:7.16}\end{equation}
in weak form. Moreover, $Z_\ell$ is monotonically increasing and 
$$
\lim_{\ell\to\infty}Z_\ell(x)
=\left(\Gamma*\left[\alpha\left(\tilde{u}^p-(u^{\kappa_*})^p\right)\right]\right)(x)
=\tilde{u}(x)-u^{\kappa_*}(x)
\quad\mbox{for a.a.~$x\in{\mathbb R}^N$}.
$$
On the other hand, by Lemma~\ref{Lemma:7.3} and \eqref{eq:7.16} we have
$$
\int_{{\mathbb R}^N}\alpha p(u^{\kappa_*})^{p-1}\varphi^{\kappa_*}Z_\ell\,dx
=\int_{{\mathbb R}^N}\nabla\varphi^{\kappa_*}\cdot\nabla Z_\ell\,dx
=\int_{{\mathbb R}^N}\alpha((u^{\kappa_*}+z_\ell)^p-(u^{\kappa_*})^p)\varphi^{\kappa_*}\,dx.
$$
Letting $\ell\to\infty$, we see that
\begin{equation}
\label{eq:7.17}
\int_{{\mathbb R}^N}
\alpha p(u^{\kappa_*})^{p-1}(\tilde{u}-u^{\kappa_*})\varphi^{\kappa_*}\,dx
=\int_{{\mathbb R}^N}\alpha\left(\tilde{u}^p-(u^{\kappa_*})^p\right)\varphi^{\kappa_*}\,dx.
\end{equation}
Let $x_0\in{\mathbb R}^N$ be such that $\tilde{u}(x_0)<\infty$. 
It follows from $\varphi^{\kappa_*}\in L^\infty_{0,-N+2}$ that
$$
\varphi^{\kappa_*}(x)\le C\omega_{0,-N+2}(x)\le C\Gamma(x_0-x)
$$
for a.a. $x\in\mathbb{R}^N$. 
Then we have
$$
\int_{{\mathbb R}^N} \alpha\tilde{u}^p\varphi^{\kappa_*}\,dx
\le C\int_{{\mathbb R}^N}\Gamma(x_0-x)\alpha\tilde{u}^p\,dx
\le C\tilde{u}(x_0)<\infty.
$$
This implies that the right-hand side of (\ref{eq:7.17}) is finite.
Since 
$$
t^p-s^p>ps^{p-1}(t-s)\quad\mbox{for $t$, $s>0$ with $t>s$},
$$
we deduce from \eqref{eq:7.17} that $\tilde{u}(x)=u^{\kappa_*}(x)$ for a.a.~$x\in{\mathbb R}^N$. 
Thus $u^{\kappa_*}$ is a unique solution to problem~\eqref{eq:P} with $\kappa=\kappa_*$.

It remains to show the equality $\kappa_*=\kappa^*$. 
Assume that $\kappa_*<\kappa^*$. 
For any $\kappa\in(\kappa_*,\kappa^*)$, 
by~\eqref{eq:2.10} we have $u^\kappa\ge(\kappa/\kappa_*)u^{\kappa_*}$. 
Then 
\begin{align*}
\int_{{\mathbb R}^N}\alpha p(u^\kappa)^{p-1}(\varphi^{\kappa_*})^2\,dx
 & \ge\left(\frac{\kappa}{\kappa_*}\right)^{p-1}
\int_{{\mathbb R}^N}\alpha p(u^{\kappa_*})^{p-1}(\varphi^{\kappa_*})^2\,dx\\
 & =\left(\frac{\kappa}{\kappa_*}\right)^{p-1}\int_{{\mathbb R}^N}|\nabla \varphi^{\kappa_*}|^2\,dx
>\int_{{\mathbb R}^N}|\nabla \varphi^{\kappa_*}|^2\,dx.
\end{align*}
On the other hand, it follows from Lemma~\ref{Lemma:5.4} that 
$$
\int_{{\mathbb R}^N}\alpha p(u^{\kappa})^{p-1}\psi^2\,dx
\le\int_{{\mathbb R}^N}|\nabla\psi|^2\,dx,
\quad \psi\in{\mathcal D}^{1,2},
$$
which is a contradiction. 
Thus the proof of Theorem~\ref{Theorem:1.2} is complete.
$\Box$
\vspace{5pt}
\newline
{\bf Proof of Corollaries~\ref{Corollary:1.1} and \ref{Corollary:1.2}.}
By \eqref{eq:1.6} we find $x_0\in{\mathbb R}^N$ such that 
$$
\limsup_{R\to 0} R^{-\frac{N}{r}}\left(\int_{B(x_0,R)}|\Gamma*\mu|^r\,dx\right)^{1/r}<\infty.
$$
Since $\alpha\equiv 1$ in ${\mathbb R}^N$, we can assume, without loss of generality, that $x_0=0$, 
which together with~\eqref{eq:1.6} implies that 
$\Gamma*\mu\in L^r_{0,\theta}$. 
Then Corollaries~\ref{Corollary:1.1} and \ref{Corollary:1.2} follow from Theorems~\ref{Theorem:1.1} and~\ref{Theorem:1.2}, respectively. 
$\Box$
\medskip

\noindent
{\bf Acknowledgment.}
The first author was supported 
in part by JSPS KAKENHI Grant Number JP19H05599. 
The second author was supported 
in part by FoPM, WINGS Program, the University of Tokyo.


\begin{bibdiv}
\begin{biblist}
\bib{Ba01}{article}{
   author={Bae, Soohyun},
   title={Asymptotic behavior of positive solutions of inhomogeneous
   semilinear elliptic equations},
   journal={Nonlinear Anal.},
   volume={51},
   date={2002},
   pages={1373--1403},
}
\bib{Ba02}{article}{
   author={Bae, Soohyun},
   title={Existence of positive entire solutions of semilinear elliptic
   equations},
   journal={Nonlinear Anal.},
   volume={71},
   date={2009},
   pages={e607--e615},
}
\bib{BCP}{article}{
   author={Bae, Soohyun},
   author={Chang, Tong Keun},
   author={Pahk, Dae Hyeon},
   title={Infinite multiplicity of positive entire solutions for a
   semilinear elliptic equation},
   journal={J. Differential Equations},
   volume={181},
   date={2002},
   pages={367--387},
}
\bib{BaL}{article}{
   author={Bae, Soohyun},
   author={Lee, Kijung},
   title={Positive entire stable solutions of inhomogeneous semilinear
   elliptic equations},
   journal={Nonlinear Anal.},
   volume={74},
   date={2011},
   pages={7012--7024},
}
\bib{BN}{article}{
   author={Bae, Soohyun},
   author={Ni, Wei-Ming},
   title={Existence and infinite multiplicity for an inhomogeneous
   semilinear elliptic equation on ${\mathbb R}^n$},
   journal={Math. Ann.},
   volume={320},
   date={2001},
   pages={191--210},
}
\bib{Be}{article}{
   author={Bernard, Guy},
   title={An inhomogeneous semilinear equation in entire space},
   journal={J. Differential Equations},
   volume={125},
   date={1996},
   pages={184--214},
}
\bib{BP}{article}{
   author={Bidaut-V\'{e}ron, Marie-Francoise},
   author={Pohozaev, Stanislav},
   title={Nonexistence results and estimates for some nonlinear elliptic
   problems},
   journal={J. Anal. Math.},
   volume={84},
   date={2001},
   pages={1--49},
}
\bib{CFY}{article}{
   author={Chen, Huyuan},
   author={Felmer, Patricio},
   author={Yang, Jianfu},
   title={Weak solutions of semilinear elliptic equation involving Dirac
   mass},
   journal={Ann. Inst. H. Poincar\'{e} C Anal. Non Lin\'{e}aire},
   volume={35},
   date={2018},
   pages={729--750},
}
\bib{CHZ}{article}{
   author={Chen, Huyuan},
   author={Huang, Xia},
   author={Zhou, Feng},
   title={Fast and slow decaying solutions of Lane-Emden equations involving
   nonhomogeneous potential},
   journal={Adv. Nonlinear Stud.},
   volume={20},
   date={2020},
   pages={339--359},
}
\bib{DGL}{article}{
   author={Deng, Yinbin},
   author={Guo, Yujin},
   author={Li, Yi},
   title={Existence and decay properties of positive solutions for an
   inhomogeneous semilinear elliptic equation},
   journal={Proc. Roy. Soc. Edinburgh Sect. A},
   volume={138},
   date={2008},
   pages={301--322},
}
\bib{DLY}{article}{
   author={Deng, Yinbin},
   author={Li, Yi},
   author={Yang, Fen},
   title={A note on the positive solutions of an inhomogeneous elliptic
   equation on ${\mathbb R}^n$},
   journal={J. Differential Equations},
   volume={246},
   date={2009},
   pages={670--680},
}
\bib{DY}{article}{
   author={Deng, Yinbin},
   author={Yang, Fen},
   title={Existence and asymptotic behavior of positive solutions for an
   inhomogeneous semilinear elliptic equation},
   journal={Nonlinear Anal.},
   volume={68},
   date={2008},
   pages={246--272},
}
\bib{Gidas}{article}{
   author={Gidas, Basilis},
   title={Symmetry properties and isolated singularities of positive
   solutions of nonlinear elliptic equations},
   book={
      series={Lecture Notes in Pure and Appl. Math.},
      volume={54},
      publisher={Dekker, New York},
   },
   date={1980},
   pages={255--273},
}
\bib{GT}{book}{
   author={Gilbarg, David},
   author={Trudinger, Neil S.},
   title={Elliptic partial differential equations of second order},
   series={Classics in Mathematics},
   note={Reprint of the 1998 edition},
   publisher={Springer-Verlag, Berlin},
   date={2001},
   pages={xiv+517},
}
\bib{HMP}{article}{
   author={Hirano, N.},
   author={Micheletti, Anna Maria},
   author={Pistoia, Angela},
   title={Multiple existence of solutions for a nonhomogeneous elliptic
   problem with critical exponent on ${\mathbb R}^N$},
   journal={Nonlinear Anal.},
   volume={65},
   date={2006},
   pages={501--513},
}
\bib{IOS01}{article}{
   author={Ishige, Kazuhiro},
   author={Okabe, Shinya},
   author={Sato, Tokushi},
   title={A supercritical scalar field equation with a forcing term},
   journal={J. Math. Pures Appl.},
   volume={128},
   date={2019},
   pages={183--212},
}
\bib{IOS02}{article}{
   author={Ishige, Kazuhiro},
   author={Okabe, Shinya},
   author={Sato, Tokushi},
   title={Thresholds for the existence of solutions to inhomogeneous elliptic equations with general exponential nonlinearity},
   journal={Advances in Nonlinear Analysis},
   volume={11},
   date={2022},
   pages={968--992},
}
\bib{IOS03}{article}{
   author={Ishige, Kazuhiro},
   author={Okabe, Shinya},
   author={Sato, Tokushi},
   title={Existence of non-minimal solutions to an inhomogeneous elliptic equation with supercritical nonlinearity},
   journal={preprint},
}
\bib{LG}{article}{
   author={Lai, Baishun},
   author={Ge, Zhihao},
   title={Infinite multiplicity for an inhomogeneous supercritical problem
   in entire space},
   journal={Proc. Amer. Math. Soc.},
   volume={139},
   date={2011},
   pages={4409--4418},
}
\bib{L}{article}{
   author={Lee, Tzong-Yow},
   title={Some limit theorems for super-Brownian motion and semilinear
   differential equations},
   journal={Ann. Probab.},
   volume={21},
   date={1993},
   pages={979--995},
}
\bib{MWL}{article}{
   author={Ma, Yong},
   author={Wang, Ying},
   author={Ledesma, C\'{e}sar T.},
   title={Lane-Emden equations perturbed by nonhomogeneous potential in the
   super critical case},
   journal={Adv. Nonlinear Anal.},
   volume={11},
   date={2022},
   pages={128--140},
}
\bib{NS01}{article}{
   author={Naito, Y\={u}ki},
   author={Sato, Tokushi},
   title={Positive solutions for semilinear elliptic equations with singular forcing terms},
   journal={J. Differential Equations},
   volume={235},
   date={2007},
   pages={439--483},
}
\bib{NS02}{article}{
   author={Naito, Y\={u}ki},
   author={Sato, Tokushi},
   title={Non-homogeneous semilinear elliptic equations involving critical
   Sobolev exponent},
   journal={Ann. Mat. Pura Appl.},
   volume={191},
   date={2012},
   pages={25--51},
}
\bib{QS}{book}{
   author={Quittner, Pavol},
   author={Souplet, Philippe},
   title={Superlinear parabolic problems},
   series={Birkh\"{a}user Advanced Texts: Basler Lehrb\"{u}cher. [Birkh\"{a}user
   Advanced Texts: Basel Textbooks]},
   publisher={Birkh\"{a}user/Springer, Cham},
   date={2019},
   pages={xvi+725},
}
\end{biblist}
\end{bibdiv}
\end{document}